\documentclass[conf]{new-aiaa}
\usepackage[utf8]{inputenc}

\usepackage{graphicx}
\usepackage{amsmath}
\usepackage[version=4]{mhchem}
\usepackage{siunitx}
\usepackage{longtable,tabularx}
\usepackage{graphicx}
\graphicspath{{Figures/}}
\usepackage{dcolumn}
\usepackage{bm}

\usepackage{float}
\usepackage{caption}
\usepackage{subcaption}

\usepackage{algorithm}
\usepackage{algorithmic}

\usepackage{soul}
\usepackage{amsmath,amsfonts}
\usepackage{graphics}
\usepackage{color}
\usepackage{xcolor}
\definecolor{rev}{rgb}{0,0,1}

\usepackage{enumitem}

\usepackage{comment}
\usepackage{hyperref}
\usepackage{lipsum}
\usepackage{tabularx}
\usepackage{stackengine}

\title{Sketching Methods for Dynamic Mode Decomposition in Spherical Shallow Water Equations}

\author{
Shady E. Ahmed\footnote{PhD Candidate, School of Mechanical and Aerospace Engineering, Oklahoma State University. shady.ahmed@okstate.edu.}, 
Omer San \footnote{Associate Professor, School of Mechanical and Aerospace Engineering, Oklahoma State University. Member AIAA. osan@okstate.edu.}
}\affil{School of Mechanical and Aerospace Engineering, Oklahoma State University, Stillwater, OK 74078, USA.}

\author{Diana A. Bistrian \footnote{Associate Professor, Department of Electrical Engineering and Industrial Informatics, Politehnica University of Timisoara. diana.bistrian@fih.upt.ro.}}
\affil{Department of Electrical Engineering and Industrial Informatics, Politehnica University of Timisoara, Hunedoara, Romania.}

\author{Ionel M. Navon \footnote{Professor, Department of Scientific Computing, Florida State University. inavon@fsu.edu.}}
\affil{Department of Scientific Computing, Florida State University, Tallahassee, Florida 32306 USA.}

\begin{document}

\maketitle

\begin{abstract}
Dynamic mode decomposition (DMD) is an emerging methodology that has recently attracted computational scientists working on nonintrusive reduced order modeling. One of the major strengths that DMD possesses is having ground theoretical roots from the Koopman approximation theory. Indeed, DMD may be viewed as the data-driven realization of the famous Koopman operator. Nonetheless, the stable implementation of DMD incurs computing the singular value decomposition of the input data matrix. This, in turn, makes the process computationally demanding for high dimensional systems. In order to alleviate this burden, we develop a framework based on sketching methods, wherein a sketch of a matrix is simply another matrix which is significantly smaller, but still sufficiently approximates the original system. Such sketching or embedding is performed by applying random transformations, with certain properties, on the input matrix to yield a compressed version of the initial system. Hence, many of the expensive computations can be carried out on the smaller matrix, thereby accelerating the solution of the original problem. We conduct numerical experiments conducted using the spherical shallow water equations as a prototypical model in the context of geophysical flows. The performance of several sketching approaches is evaluated for capturing the range and co-range of the data matrix. The proposed sketching-based framework can accelerate various portions of the DMD algorithm, compared to classical methods that operate directly on the raw input data. This eventually leads to substantial computational gains that are vital for digital twinning of high dimensional systems.
\end{abstract}

\section{Introduction} \label{sec:intro}
With the advent of relatively inexpensive and powerful sensors and enhanced data storage capabilities, data sets from various processes are ubiquitous in all fields of engineering. However, processing the data and extracting meaningful information is the bottleneck of the entire data handling pipeline. Despite the huge progress in high performance computing, computational power remains a limiting factor with respect to the amount of information that is generated and processed on a daily basis. To mitigate the load on computational analysis, approaches to extract useful information from the entire data pool are sought. The idea is that not every piece of information is needed to analyse a process. Beneath all high-dimensional data lies considerably lower dimensional patterns which govern most of the dynamics. Modal decomposition techniques have been developed to identify such dominant spatio-temporal modes which dominate the evolution of the system \cite{taira2017modal, gueniat2015dynamic, puzyrev2019pyrom}. This introduces the concept of reduced order models (ROMs) where the entire system can be well represented by relatively compact surrogates, opening ways for computationally inexpensive analysis of the prominent dynamics \cite{quarteroni2015reduced, rowley2017model}.  

Among the present methodologies that we have at our disposal for ROMs, dynamic mode decomposition (DMD) has rapidly gained recognition during the last few years \cite{schmid2010dynamic,rowley2009spectral,tissot2014model}. DMD shares some roots with the Koopman theory and can be viewed as a data-driven approximation of the Koopman operator spectrum. It has some interesting advantages, especially if the underlying dynamics is quasi-periodic and well characterized by fast decaying singular values. For one, it is completely data driven and does not require any prior knowledge of the dynamics of the concerned system. This makes it simple to use with quite complicated dynamical systems. It can be applied to both experimental as well as numerical data with success. Furthermore, it is theoretically sound and different analyses can be applied to DMD as it is based on the Koopman operator. 


Since its introduction in fluid dynamics community \cite{schmid2010dynamic,rowley2009spectral}, DMD has been one of the most widely used reduced order modeling strategies in many other disciplines as well \cite{erichson2019compressed,natsume2020application}. With increasing popularity, various strategies have been developed to overcome many of the early DMD shortcomings. Studies to improve its performance via pre-processing and post-processing have been conducted like ensemble-averaging methods \cite{sarmast2014mutual}, mean subtraction \cite{hirsh2020centering}, de-biasing algorithms \cite{hemati2017biasing}, sparsity inducing approaches \cite{jovanovic2014sparsity}, online update techniques \cite{hemati2014dynamic,zhang2019online}, and improved least square methods \cite{brunton2014compressive,kutz2016multiresolution}. Moreover, linear and nonlinear aspects of DMD have been discussed in \cite{alekseev2016linear}. An extended dynamic mode decomposition approach has been proposed to better approximate the Koopaman operator tuples (i.e., eigenvalues, eigenfunctions, and modes) \cite{williams2015data}. In addition, efforts have been made to make the DMD algorithm robust against noisy data \cite{schmid2010dynamic,scherl2020robust}. 

Another important aspect of reduced order modeling by modal decomposition techniques such as DMD involves effective selection of the modes. In proper orthogonal decomposition (POD) based model reduction approaches \cite{sirovich1987turbulence,amsallem2008interpolation,holmes2012turbulence}, modes are automatically sorted based on the eigenvalues of the auto-correlation matrix of snapshot data. Those eigenvalues intrinsically represent the amount of total system's energy captured by the individual modes. In contrast to the POD based approaches, DMD does not inherently rank the underlying modes and the selection criteria is not unique nor straightforward \cite{chen2012variants,jovanovic2014sparsity,sayadi2015parametrized}. Initial implementations \cite{rowley2009spectral,wan2015dynamic} utilized the norm of modes as a ranking and sorting criteria. \citet{sayadi2014reduced} used the projection of the first snapshot onto the modes to classify them. However, this method does not perform well for fast-decaying modes. \citet{jovanovic2014sparsity} provided the computation of optimal DMD amplitudes based on the minimization of the difference between snapshot data and DMD reconstruction over the total time window. Similarly, an optimized DMD algorithm was proposed to compute arbitrary number of DMD modes to improve DMD accuracy when a few modes are sought \cite{chen2012variants}. Ideas from data assimilation, including variational and sequential approaches, have been also adopted to build physically-sound DMD-based ROMs \cite{tissot2014model,tissot20154d,nonomura2018dynamic,nonomura2019extended}.

The straightforward implementation approach proposed by \citet{rowley2010reduced} involves the computation of a companion matrix that helps to construct, in a least squares sense, the final data vector as a linear combination of all previous data vectors. However, \citet{schmid2010dynamic} showed that this version may be ill-conditioned in practice and an alternative algorithm based on the singular value decomposition (SVD) of the snapshot data matrix is recommended. Nonetheless, the brute-force computation of the SVD for high dimensional systems becomes computationally prohibitive. In this study, we explore the applicability sketching methods for the efficient computation of DMD basis and spectrum. These methods rely upon informed projections and aim at constructing a smaller (memory-wise) matrix, called the ``sketch'' while preserving important properties of the original data \cite{mahoney2011randomized,woodruff2014sketching}. The sketch may represent any combination of row space, column space or the space generated by the intersection of rows and columns (core space). There have been several studies that aimed to utilize random projections to derive low-rank randomized DMD \cite{erichson2016randomized,bistrian2017randomized,bistrian2018efficiency,erichson2019randomized,alla2019randomized}. The major contribution of this study is to build on and extend these previous efforts. In particular, we explore the applicability of sketching algorithms in the low-rank DMD computations for spherical shallow water equations data, representing a relatively simple, but very representative model for geophysical flows. The spherical SWEs formulation is often considered a first step in developing general circulation models in large scales. We also explore a sketching algorithm that aims at capturing the range and corange of the snapshot data matrix with free parameters that can be tuned based on the given memory constraints. In addition, we investigate the effect of different sorting criteria onto the accuracy of the low-rank DMD reconstruction. We further show that sketching-based DMD can directly yield a near-target-rank DMD and mitigate the need for special mode-selection criterion.

The rest of the paper is structured as follows. In Section~\ref{sec:math}, we present the underlying governing equations of shallow water system on spherical coordinates and the corresponding initial and boundary conditions. Numerical schemes used to solve the governing equations for data generation as well metrics for evaluating results are summarized in Section~\ref{sec:numerics}. We then detail the dynamic mode decomposition formulations and the considered sketching approaches in Section~\ref{sec:DMD}. Results are provided with corresponding discussions in Section~\ref{sec:res}, while concluding remarks are drawn in Section~\ref{sec:conc}.

\section{Mathematical Model } \label{sec:math}
The shallow water equations (SWEs) are the mass and momentum balance equations that constitute a specialized case of the Navier-Stokes equations (NSEs). The SWEs describe the flow field of a free fluid surface in cases where the horizontal length scale dominates over the vertical length scale; implying that the horizontal velocity field is approximately invariant in the vertical direction. Thus, the variation of the vertical component vanishes in the SWEs. Mathematically, SWEs are obtained by the integrated average of the NSEs across the vertical length and substituting the pressure term with the depth of the fluid column through the hydrostatic approximation. The Coriolis term is included to account for the forces due to Earth's rotation in cases of geophysical flows. Interestingly, in many atmospheric and oceanic flows of practical interest, such assumptions hold and hence, they can be adequately modeled by these equations \cite{williamson1992standard,thuburn2000numerical}. Therefore, SWEs form a good test bed for reduced order modeling algorithms in the context of geophysical flows. The SWEs for the atmosphere on Earth, using spherical coordinates, can be written as follows:

\begin{equation} \label{eq:SWE1}
    \frac{\partial h}{\partial t} + \frac{1}{\rho \cos{\theta}} \frac{\partial h u_{\phi}}{\partial \phi} + \frac{1}{\rho} \frac{\partial h u_{\theta}}{\partial \theta} = \frac{h u_{\theta}}{\zeta} \tan{\theta},
\end{equation}

\begin{equation} \label{eq:SWE2}
    \frac{\partial h u_\phi}{\partial t} + \frac{1}{\rho \cos{\theta}} \frac{\partial}{\partial \phi} \left( h u_{\phi}^2 + \frac{1}{2} g h^2 \right) + \frac{1}{\rho}\frac{\partial \left( h u_{\phi} u_{\theta} \right)}{\partial \theta} = Fhu_{\theta} - \frac{gh}{\zeta \cos{\theta}} \frac{\partial H}{\partial \phi} + \frac{h u_{\phi} u_{\theta}}{\zeta} \tan{\theta},
\end{equation}

\begin{equation} \label{eq:SWE3}
    \frac{\partial h u_\theta}{\partial t} + \frac{1}{\rho \cos{\theta}} \frac{\partial \left( h u_{\phi} u_{\theta} \right)}{\partial \phi} + \frac{1}{\rho}\frac{\partial}{\partial \theta} \left( h u_{\theta}^2 + \frac{1}{2} gh^2 \right) = -Fhu_{\phi} - \frac{gh}{\zeta} \frac{\partial H}{\partial \theta} + \frac{h u_{\phi}^2}{\zeta} \tan{\theta}.
\end{equation}

In the above system of equations, $\zeta=\rho+H$ where $\rho$ denotes the Earth's radius ($\rho = 6.4 \times 10^{6}$ m) and $H$ is the height of the bottom topography; $g$ is the acceleration due to gravity ($g=9.8$ m/s$^2$); $h$ is the depth of the water surface; $u_\theta$ and $u_\phi$ are the respective velocities in the longitudinal and latitudinal directions, while $\phi$ and $\theta$ are the longitudes and latitudes, respectively. 

We consider a spherical domain given by longitudes $\Phi = [0,360^{\circ}]$ and latitudes $\Theta = [-79.5^{\circ},79.5^{\circ}]$ for the purpose of our study. For simplicity, the bottom surface is taken to be flat, i.e, $H(\phi,\theta)=0 \; \forall \; \phi \in \Phi$ and $\theta \in \Theta$. Initial condition for $h$ is given by:
\begin{equation} \label{eq:IC1}
    h(\theta) = 10000 - 60 \cos{\left( 4\pi\theta \right)}e^{-\theta^2};  \;\;\; \forall \phi \in \Phi.
\end{equation}
In order to trigger shear layer instability, we add random disturbance to the initial height as follows:
\begin{equation} \label{eq:IC1.1}
    h(\theta) = 10000 - 60 \cos{\left( 4\pi\theta \right)}e^{-\theta^2} + \kappa \frac{\Delta \theta (N_{\theta}-1)}{\pi} |F| \times 10^{4} \cos{\theta}, 
\end{equation}
where $\kappa$ is a uniformly distributed random number between 0 and 1, $\Delta \theta$ is the spatial resolution in the latitudinal direction (one degree in this study), $N_{\theta}$ is the number of latitudinal grid points (i.e., $160$ points), and $F=2f\sin{\theta}$ is the Coriolis parameter with $f=\frac{2\pi}{24 \times 3600}$ rad/s being the Earth's rotation rate.

For initial velocities, geostrophic wind conditions are assumed to prevail and accordingly the initial $u^0_{\theta}$ and $u^0_{\phi}$ are given as:
\begin{equation} 
    u_{\phi}(\phi,\theta) = - \frac{g}{\zeta(F-\delta)} \frac{\partial h}{\partial \theta},\quad
    u_{\theta}(\phi,\theta) = \frac{g}{\zeta(F-\delta) \cos{\theta}} \frac{\partial h}{\partial \phi}, \qquad \forall \phi \in \Phi,  \ \forall \theta \in \Theta \label{eq:IC2.2}
\end{equation}
where $\delta$ is a small constant introduced to prevent the value the denominator from becoming exactly $0$ at the equator ($\theta=0^{\circ}$). Here, the first order spatial derivative of $h$ is numerically computed from the initial condition given by Eq.~\ref{eq:IC1} using central difference scheme. We use periodic boundary conditions in the longitudinal ($\phi$) direction and slip boundary conditions in the latitudinal ($\theta$) direction as given below for the boundary grid points:
\begin{gather} 
    h(\phi,\theta,t) = h(\phi+\phi_L,\theta,t), \label{eq:BC1.1} \\
    u_{\phi}(\phi,\theta,t) = u_{\phi}(\phi+\phi_L,\theta,t), \label{eq:BC1.2} \\
    u_{\theta}(\phi,\theta,t) = u_{\theta}(\phi+\phi_L,\theta,t), \label{eq:BC1.3}
\end{gather}

\begin{gather} 
    \frac{\partial h}{\partial \theta}\bigg|_{(\phi,\theta_0,t)} = 0, \;\;\;\;\;\;    
    \frac{\partial h}{\partial \theta}\bigg|_{(\phi,\theta_L,t)} = 0, \label{eq:BC2.1} \\
    \frac{\partial u_{\phi}}{\partial \theta}\bigg|_{(\phi,\theta_0,t)} = 0, \;\;\;\;\;\;    
    \frac{\partial u_{\phi}}{\partial \theta}\bigg|_{(\phi,\theta_L,t)} = 0, \label{eq:BC2.2} \\
    u_{\theta}(\phi,\theta_0,t) = 0, \;\;\;\;\;\;            
    u_{\theta}(\phi,\theta_L,t) = 0. \label{eq:BC2.3}
\end{gather}
where $\phi_L=360^{\circ}$, $\theta_0=-79.5^{\circ}$ and $\theta_L=79.5^{\circ}$ (corresponding to the domain boundaries). 

\section{Numerical Methods} \label{sec:numerics}

We use the Lax-Wendroff method for solving the set of partial differential equations of the SWEs on the sphere. It is a second-order accurate scheme in both time and space. Consider a generic variable $q$, which has a conservation equation defined as:

\begin{equation} \label{eq:gen_con}
    \frac{\partial \mathbf{q}}{\partial t} + \frac{1}{\rho\cos{\theta}} \frac{\partial \mathbf{f}(\mathbf{q})}{\partial \phi} + \frac{1}{\rho} \frac{\partial \mathbf{g}(\mathbf{q})}{\partial \theta} = \mathbf{Q}(\mathbf{q}). 
\end{equation}
In our case, $\mathbf{q} = [h,hu_{\phi},hu_{\theta}]^T$ refers to the vector of conservative variables. Comparing Eq.~\ref{eq:gen_con} with each of Eq.~\ref{eq:SWE1}, Eq.~\ref{eq:SWE2} and Eq.~\ref{eq:SWE3}, we obtain the fluxes and the source vector as follows:
\begin{gather}
    \mathbf{f}(\mathbf{q}) = \left[h u_{\phi} \; , \;  hu^2_{\phi}+\frac{1}{2}gh^2 \; , \; hu_{\phi}u_{\theta} \right]^T, \\
    \mathbf{g}(\mathbf{q}) = \left[hu_{\theta} \; , \; hu_{\phi}u_{\theta} \; , \; hu^2_{\theta}+\frac{1}{2}gh^2 \right]^T, \\
    \mathbf{Q}(\mathbf{q}) = h\left[ \alpha \; , \; \beta\; , \; \gamma \right]^T.
\end{gather} 
where
\begin{gather}
\alpha = \frac{u_{\theta}\tan{\theta}}{\zeta}, \\
\beta = Fu_{\theta}-\frac{g}{\zeta\cos{\theta}} \frac{\partial H}{\partial \phi} + \frac{u_{\phi}u_{\theta}}{\zeta}\tan{\theta}, \\
\gamma =  -Fu_{\phi}-\frac{g}{\zeta} \frac{\partial H}{\partial \theta} + \frac{u_{\phi}^2}{\zeta}\tan{\theta} .
\end{gather} 
Now, the first step in the Lax-Wendroff scheme is to estimate the value of $\mathbf{q}$ at the midpoints in space and time as:
\begin{gather}
    \mathbf{q}^{(k+1/2)}_{i+1/2,j} = \frac{\mathbf{q}^{(k)}_{i,j} + \mathbf{q}^{(k)}_{i+1,j}}{2} - \frac{\Delta t}{2} \frac{1}{\rho\cos{\theta}} \left(\frac{\mathbf{f}^{(k)}_{i+1,j}-\mathbf{f}^{(k)}_{i,j}}{\Delta \phi} \right), \\
    \mathbf{q}^{(k+1/2)}_{i,j+1/2} = \frac{\mathbf{q}^{(k)}_{i,j} + \mathbf{q}^{(k)}_{i,j+1}}{2} - \frac{\Delta t}{2} \frac{1}{\rho} \left(\frac{\mathbf{g}^{(k)}_{i,j+1}-\mathbf{g}^{(k)}_{i,j}}{\Delta \theta} \right),
\end{gather}
where the superscript $(k)$ denotes the time index corresponding to discrete time $t_k$ while the subscripts define the spatial index with respect to $\phi$ and $\theta$ directions, respectively. $\Delta t$ is the time step size while $\Delta \phi$ and $\Delta \theta$ are the corresponding spatial resolutions of the spherical grid. The value of $\mathbf{q}$ at next time step is given in the second step as follows:
\begin{equation}
    \Tilde{\mathbf{q}}^{(k+1)}_{i,j} = \mathbf{q}^{(k)}_{i,j} + \Delta t \left[ - \left( \frac{1}{\rho\cos{\theta}}\right) \frac{\mathbf{f}^{(k+1/2)}_{i+1/2,j} - \mathbf{f}^{(k+1/2)}_{i-1/2,j}}{\Delta \phi} - \left( \frac{1}{\rho} \right)  \frac{\mathbf{g}^{(k+1/2)}_{i,j+1/2} - \mathbf{g}^{(k+1/2)}_{i,j-1/2}}{\Delta \theta} \right].
\end{equation}

As a last step, before we move to the next iteration, the primitive quantities are updated by adding the source terms as follows:
\begin{gather}
h^{(k+1)}_{i,j} = \tilde{h}^{(k+1)}_{i,j} + \frac{\Delta t}{2} \alpha^{(k)}_{i,j} \left( \tilde{h}^{(k+1)}_{i,j} + h^{(k)}_{i,j}\right), \\
 {u_{\phi}}^{(k+1)}_{i,j} = \frac{1}{h^{(k+1)}_{i,j}} \left[ (\tilde{h} \tilde{u}_{\phi})^{(k+1)}_{i,j}  + \frac{\Delta t}{2} \beta^{(k)}_{i,j} \left( h^{(k+1)}_{i,j} + h^{(k)}_{i,j}\right) \right],  \\
  {u_{\theta}}^{(k+1)}_{i,j} = \frac{1}{h^{(k+1)}_{i,j}} \left[ (\tilde{h} \tilde{u}_{\theta})^{(k+1)}_{i,j}  + \frac{\Delta t}{2} \gamma^{(k)}_{i,j} \left( h^{(k+1)}_{i,j} + h^{(k)}_{i,j}\right) \right].  
\end{gather} 

Although more sophisticated numerical methods that address inherent numerical difficulties are available for solving spherical SWE \cite{neta1997analysis}, in this work, without considering poles (i.e., $\Theta = [-79.5^{\circ},79.5^{\circ}]$), we simply deploy the described Lax-Wendroff scheme to numerically solve the SWEs on a sphere to obtain the solution of the surface flow dynamics \cite{connolly2017}. This method is applied to each of the three equations of spherical SWEs simultaneously to obtain the height ($h$) as well as the velocity ($u_{\phi}$ and $u_{\theta}$) fields. For the forward simulations, we use spatial resolution $\Delta \theta = 1^{\circ}$ and $\Delta \phi = 1^{\circ}$ and a time step of $\Delta t = 30$ sec. The total time duration for the simulation is set to $6$ days and the snapshots of the field data are stored at every $15$ minutes. Moreover, we truncate the snapshots corresponding to the first $3$ days from our input data matrix to let the system pass the initial transition period. Thus, we get data set in the form of $289$ snapshots arrays of dimensions $360 \times 160$ (corresponding to the spatial resolution) for each flow field, where each snapshot corresponds to a particular time. In the present study, we consider the vorticity field data, which is defined as the curl of the velocity vector, i.e.,
\begin{equation} \label{eq:vorticity1}
    \boldsymbol{\omega} = \nabla \times \mathbf{u}.
\end{equation}
Eq.~\ref{eq:vorticity1} in the spherical coordinate can be simplified for the 2D case as:
\begin{equation} \label{eq:vorticity2}
    \omega = \frac{1}{\rho \cos{\theta}} \left[ \frac{\partial u_{\theta}}{ \partial \phi} - \frac{\partial}{\partial \theta} \left( u_{\phi} \cos{\theta} \right) \right].
\end{equation}
Eq.~\ref{eq:vorticity2} is discretized using central difference scheme at the central grid points. For the boundary grid points, we similarly use periodic boundary condition in the longitudinal direction and slip condition in the latitudinal direction.

\section{Dynamic Mode Decomposition} \label{sec:DMD}
Dynamic mode decomposition (DMD) is one of the most popular methods for data-based reduced order modeling. The core DMD framework by itself is completely data-driven and does not require prior knowledge of the model dynamics. Let $\mathbf{x} \in \mathbb{R}^n$ (where $n \gg 1$) be the state of a system that evolves in time through some specific dynamics as $\dfrac{d \mathbf{x}}{dt} = f(\mathbf{x},t)$, where $f(\cdot)$ can be a nonlinear function. The objective of DMD is to identify the leading DMD eigenvalues and the corresponding eigenvectors of the best fit operator $\mathbf{A}$ which would evolve the system linearly in time as:
\begin{equation} \label{eq:cont_DMD_op}
    \frac{d \mathbf{x}}{dt} \approx \mathbf{A} \mathbf{x}.
\end{equation}

By estimating the operator $\mathbf{A}$, spatially coherent DMD modes can be computed, each of which is associated with a growth/decay rate and a frequency that define its time dynamics. Dimensionality reduction is obtained by choosing a few significant modes while rejecting the others. DMD was formally introduced to the fluid community by \citet{schmid2011applications}. In problems of fluid dynamics, the system dimension $n$ can be very large to account for the relevant scales in the flow field, and frameworks such as DMD are of high interest. In our case, we look at SWEs on the sphere as our governing dynamical system (described in Section \ref{sec:math}). 

In practice, one would have data of the system state $\mathbf{X} = [ \mathbf{x}^{(1)} , \mathbf{x}^{(2)} , \dots , \mathbf{x}^{(m)}]  \in \mathbb{R}^{n \times m}$ collected at various time instants $t_k = [t_1,t_2,\dots,t_m] \in \mathbb{R}^m$. So, here $n$ is the number of degrees of freedom in the system while $m$ is the number of discrete time steps for which data is available (i.e., number of collected snapshots). This data can be either actual field data from physical experiments or even synthetic data generated from high fidelity numerical computations. In the present study, data is obtained from the solution of the governing equations as described in Section \ref{sec:numerics}, which constitute the full order model (FOM) of the system of interest. We focus on the vorticity field on the spherical domain $[\Phi \times \Theta]$. Each vorticity snapshot matrix of dimension $n=N_{\phi} \times N_{\theta}$ (i.e., $360 \times 160$), corresponding to a particular time $t_k$, is rearranged in a column vector $\mathbf{x}^{(k)}\in \mathbb{R}^{360 \cdot 160}$. So, we have a dynamical system where the system state $\mathbf{x}^{(k)} = \mathbf{x}(t_k) \in \mathbb{R}^{57600} (1 \leqslant k \leqslant m)$ denotes the vorticity field on a spherical domain $[\Phi \times \Theta]$ at a particular time $t_k$. Thus, the full data matrix $\mathbf{x} \in \mathbb{R}^{n \times m}$ ($n=57600$ for our case) is formed as:
\begin{equation} \label{eq:full_data_mat}
 \mathbf{X} = \begin{bmatrix}
 | & | &   & | \\
 \mathbf{x}^{(1)} & \mathbf{x}^{(2)} & \dots & \mathbf{x}^{(m)} \\ 
 | & | &  & | \\
 \end{bmatrix}.
\end{equation}
The primary objective of DMD is to extract the eigenvectors and eigenvalues of the matrix $\hat{\mathbf{A}}$ such that:
\begin{equation} \label{eq:dist_DMD_op}
    \mathbf{x}^{(k+1)} \approx \hat{\mathbf{A}} \mathbf{x}^{(k)},
\end{equation}
where Eq.~\ref{eq:dist_DMD_op} is the time discretized form of Eq.~\ref{eq:cont_DMD_op}. In essence, both are similar in that they denote the evolution of the system states as a linear operation. Next, the full data set is split into two matrices $\mathbf{X}_1 \in \mathbb{R}^{n \times (m-1)}, \mathbf{X}_2 \in \mathbb{R}^{n \times (m-1)}$ as defined below:
\begin{equation} \label{eq:split_data_mat}
 \mathbf{X}_1 = \begin{bmatrix}
 | & | &   & | \\
 \mathbf{x}^{(1)} & \mathbf{x}^{(2)} & \dots & \mathbf{x}^{(m-1)} \\ 
 | & | &  & | \\
 \end{bmatrix}, \qquad
 \mathbf{X}_2 = \begin{bmatrix}
 | & | &   & | \\
 \mathbf{x}^{(2)} & \mathbf{x}^{(3)} & \dots & \mathbf{x}^{(m)} \\ 
 | & | &  & | \\
 \end{bmatrix}.
\end{equation}
Thus, Eq.~\ref{eq:dist_DMD_op} ca be rewritten in matrix format as follows:
\begin{equation} \label{eq:DMD_op}
    \mathbf{X}_2 \approx \hat{\mathbf{A}} \mathbf{X}_1.
\end{equation}
The operator $\hat{\mathbf{A}}$, representing the best linear fit, can be computed through the least squares optimization of:
\begin{equation} \label{eq:optimization}
    \hat{\mathbf{A}} = \underset{\hat{\mathbf{A}}}{\arg\min} \|\mathbf{X}_2-\hat{\mathbf{A}}\mathbf{X}_1 \|_F \;,
\end{equation}
where $\|.\|_F$ is the Frobenius norm. In the present study, we investigate different ways to approximate the eigenvalues and eigenvectors to form DMD modes. Ideally, the full rank $\hat{\mathbf{A}}$ matrix would have $\mathbb{R}^{n \times n}$ dimension, which would require high computational power to operate when $n$ is very large (which is the case in practical fluid dynamics problems). So, the motivation here is to seek ways to replace it with lower rank approximations. This will be described in detail in Section~\ref{sub:std_DMD}

Once the DMD modes are computed, we can form a matrix $\boldsymbol{\Psi} = \{\boldsymbol{\psi}_i\}_{i=1}^m \in \mathbb{C}^{n \times m}$ where each column $\boldsymbol{\psi_i}$ represents a DMD mode. DMD inherently does not provide robust ranking of the modes unlike other algorithms such as proper orthogonal decomposition (POD). So, an additional criterion has to be incorporated to appropriately rank these modes. The next step is to truncate the DMD mode matrix $\boldsymbol{\Psi}$ to $\boldsymbol{\Psi}^{}_\mathbf{r} = \{ \boldsymbol{\psi}_i \}_{i=1}^{r} \in \mathbb{C}^{n \times r}$ where $r$ ($\ll m$) is the number of retained modes. The information of the time dynamics in each DMD basis is inferred from the corresponding eigenvalues $\boldsymbol{\Lambda} = \text{diag}\{ \lambda_i \}_{i=1}^{r} \in \mathbb{C}^{r \times r}$. Here, $\lambda_i$ are the discrete-time eigenvalues, while the continuous-time eigenvalues can be evaluated as follows:
\begin{equation}\label{eq:dis_cont_evalue}
    \alpha_i = \frac{\ln{\lambda_i} }{\Delta T},
\end{equation}
where $\Delta T$ is the time interval between two consecutive snapshots and $\boldsymbol{\alpha} = \{ \alpha_i \}_{i=1}^{r} \in \mathbb{C}^{r}$ constitutes the vector of the continuous-time eigenvalues. Next, we reconstruct the high dimensional dynamics from the lower order subspace as:
\begin{equation} \label{eq:recon_dis}
    \mathbf{x}^{(k)}_{ROM} = \sum_{i=1}^{r}\boldsymbol{\psi}_i \lambda_i^{k-1} b_i = \boldsymbol{\Psi}^{}_\mathbf{r} \boldsymbol{\Lambda}^{k-1} \mathbf{b},
\end{equation}
where $\mathbf{b} = \{ b_i \}_{i=1}^{r} \in \mathbb{C}^{r}$ is the vector of initial amplitudes of the DMD modes given as:
\begin{equation} \label{eq:amp}
     \mathbf{b} = \boldsymbol{\Psi}^{\dagger}_\mathbf{r} \mathbf{x}^{(1)},
\end{equation}
where $\boldsymbol{\Psi_{r}^{\dagger}}$ is the pseudoinverse of $\boldsymbol{\Psi}^{}_\mathbf{r}$. 
In Eq.~\ref{eq:recon_dis}, the superscript on $\mathbf{x}$ is the index for discrete time (i.e., $\mathbf{x}^{(k)} = \mathbf{x}(t_k)$), while the superscripts (without parenthesis) on $\lambda_i$ and $\boldsymbol{\Lambda}$ refer to powers. Finally, for reconstruction, Eq.~\ref{eq:recon_dis} can be written in terms of the continuous-time eigenvalues by using the relation in Eq.~\ref{eq:dis_cont_evalue} as follows:
\begin{equation} \label{eq:recon_cont}
     \mathbf{x}^{(k)}_{ROM} = \sum_{i=1}^{r} \boldsymbol{\psi}_i e^{\alpha_i t_k} b_i = \boldsymbol{\Psi}^{}_\mathbf{r} \text{diag}[e^{\boldsymbol{\alpha}t_k}] \mathbf{b}.
\end{equation}
 
In Eq.~\ref{eq:recon_cont}, the real part of the continuous-time eigenvalues is responsible for the growth/decay rate of the mode. The mode grows for a positive real part and conversely decays over time for a negative real part. On the other hand, the imaginary part determines the oscillating frequency of the mode.
 
\subsection{Deterministic DMD} \label{sub:std_DMD}
In the standard form of the DMD, the linear operator $\hat{\mathbf{A}}$ is projected onto a lower $R$-dimensional subspace to replace the full dimensional $\hat{\mathbf{A}}$ matrix. The projection of $\hat{\mathbf{A}}$ is referred to as $\Tilde{\mathbf{A}}$. The first step is to take the singular value decomposition (SVD) of the matrix $\mathbf{X}_1$ as $\mathbf{X}_1 = \mathbf{U}\boldsymbol{\Sigma}\mathbf{V}^*$ (where $\mathbf{V}^*$ denotes the complex conjugate transpose of matrix $\mathbf{V}$). Compact SVD  can be computed such that $\mathbf{X}_1 \approx \mathbf{U}^{}_\mathbf{R} \boldsymbol{\Sigma}^{}_\mathbf{R} \mathbf{V}^*_\mathbf{R}$, where $\mathbf{U}^{}_\mathbf{R} \in \mathbb{R}^{R \times R}$ and $\mathbf{V}^{}_\mathbf{R} \in \mathbb{R}^{R \times R}$ denote the matrix of the first $R$ columns of $\mathbf{U}$ and $\mathbf{V}$ respectively, while $\boldsymbol{\Sigma}^{}_\mathbf{R} \in \mathbb{R}^{R \times R}$ is the first $R \times R$ dimensional sub-block of $\boldsymbol{\Sigma}$, with $R$ being the rank of $\boldsymbol{\Sigma}$.

The projection of $\hat{\mathbf{A}}$ onto the $R$-dimensional space is taken as:
\begin{equation} \label{eq:Atilde}
    \Tilde{\mathbf{A}} = \mathbf{U}^*_\mathbf{R} \hat{\mathbf{A}} \mathbf{U}^{}_\mathbf{R}.
\end{equation}
Using Eq.~\ref{eq:Atilde}, the optimization problem in Eq.~\ref{eq:optimization} becomes $\Tilde{\mathbf{A}} = \underset{\Tilde{\mathbf{A}}}{\arg\min} \| \mathbf{X}_2-\mathbf{U^{}_R} \Tilde{\mathbf{A}} \boldsymbol{\Sigma}^{}_\mathbf{R} \mathbf{V}^*_\mathbf{R} \|_F$, which gives $\Tilde{\mathbf{A}} = \mathbf{U^*_R} \mathbf{X}_2 \mathbf{V^{}_R} \boldsymbol{\Sigma}^{-1}_{\mathbf{R}} \in \mathbb{R}^{R \times R}$ . The eigenvalues of $\Tilde{\mathbf{A}}$ represent the Ritz values (approximate eigenvalues) of $\hat{\mathbf{A}}$. 
These eigenvalues (also called DMD eigenvalues) are computed as $\Tilde{\mathbf{A}} \mathbf{W} = \mathbf{W} \boldsymbol{\Lambda}$. It should be noted that the eigen decomposition of $\Tilde{\mathbf{A}}$ results in complex eigenvalues and eigenvectors. Here, each columns of $\mathbf{W} \in \mathbb{C}^{R \times R}$ are the eigenvectors while $\boldsymbol{\Lambda} = \text{diag} (\{ \lambda_i \}_{i=1}^R) \in \mathbb{C}^{R \times R}$ is the diagonal matrix containing the corresponding eigenvalues. The DMD modes can be evaluated as follows:
\begin{equation}
    \boldsymbol{\Psi} = \mathbf{U}_\mathbf{R}\mathbf{W},
\end{equation}
while the so-called ``exact'' DMD modes can be written as \cite{kutz2016dynamic,tu2014dynamic}:
\begin{equation}
    \boldsymbol{\Psi} = \mathbf{X}_2 \mathbf{V}^{}_\mathbf{R} \boldsymbol{\Sigma}^{-1}_\mathbf{R} \mathbf{W}.
\end{equation}
A suitable reduced order approximation $\boldsymbol{\Psi}^{}_\mathbf{r}$ can be obtained by retaining $r$ columns of $\boldsymbol{\Psi}$. From here, the continuous-time eigenvalues can be easily calculated from Eq.~\ref{eq:dis_cont_evalue} and the estimated higher order dynamics from the lower order approximation can be found from either Eq.~\ref{eq:recon_dis} or Eq.~\ref{eq:recon_cont}. Algorithm~\ref{alg:std_DMD} shows the standard SVD-based DMD implementation.

\begin{algorithm}[H]
\caption{Deterministic Dynamic Mode Decomposition}
\label{alg:std_DMD}
\begin{algorithmic}[1]
\STATE The matrix $\mathbf{X}$ is split into two matrices $\mathbf{X}_1 = \{ \mathbf{x}^{(1)}, \mathbf{x}^{(2)}, \dots, \mathbf{x}^{(m-1)} \}$ and $\mathbf{X}_2 = \{ \mathbf{x}^{(2)}, \mathbf{x}^{(3)}, \dots, \mathbf{x}^{(m)} \}$. \\
\vspace{10pt}
\STATE Perform SVD on $\mathbf{X}_1$  \vspace{-10pt}
\begin{align*}
    \mathbf{U} \boldsymbol{\Sigma} \mathbf{V}^* = \text{svd} {\left( \mathbf{X}_1 \right)}
\end{align*}

\STATE Rank truncation [to reduce noise] \vspace{-10pt}
\begin{align*}
    \mathbf{U}_\mathbf{R} &= \mathbf{U}(:,1:R) \\
    \mathbf{V}_\mathbf{R} &= \mathbf{V}(:,1:R) \\
    \boldsymbol{\Sigma}_\mathbf{R} &= \boldsymbol{\Sigma}(1:R,1:R)
\end{align*}

\STATE Low-rank dynamics \vspace{-10pt}
\begin{align*}
    \Tilde{\mathbf{A}} = \mathbf{U^*_R} \mathbf{X}_2 \mathbf{V}^{}_\mathbf{R} \boldsymbol{\Sigma}^{-1}_\mathbf{R}
\end{align*}

\STATE Eigenvalue decomposition \vspace{-10pt}
\begin{align*}
    \left[\mathbf{W}, \boldsymbol{\Lambda} \right] = \text{eig} {( \Tilde{\mathbf{A}} )}
\end{align*}

\STATE Compute DMD modes and spectrum \vspace{-10pt}
\begin{align*}
\boldsymbol{\Psi} &= \mathbf{U}_\mathbf{R} \mathbf{W} \quad \text{or} \quad \boldsymbol{\Psi} = \mathbf{X}_2 \mathbf{V}^{}_\mathbf{R} \boldsymbol{\Sigma}^{-1}_\mathbf{R} \mathbf{W} \\
\lambda_i &= \{ \text{diag} (\boldsymbol{\Lambda}) \} \\
\alpha_i &= \text{ln} {(\lambda_i)}/\Delta T
\end{align*}
\end{algorithmic}
\end{algorithm} 

\subsection{Rank truncation and mode selection}
Mode selection forms a crucial part of ROMs since the entire concept of ROMs is premised upon the concept of predicting the system dynamics based on the projection of the full order dynamics onto a comparatively small number of bases (called modes). Although the DMD results in modes each of which is associated with a unique oscillation frequency and growth/decay rate, it does not possess an inherent sorting criteria. Therefore, user-specific rules have to be set to select the most influential modes, and here we discuss some of the approaches to do so.

\subsubsection{Early truncation}
One can simply utilize the fact that the DMD computation involves a singular value decomposition, resulting in hierarchically arranged bases for the snapshot data matrix $\mathbf{X}_1$. Therefore, an early truncation can be performed by using an $r$-rank approximation of $X_1$ by setting $R=r$ in Section~\ref{sub:std_DMD}. This early truncation eventually yields exactly $r$ DMD modes, with no need for further sorting or selection mechanisms. Nonetheless, it has been shown in previous studies (e.g., \cite{ahmed2020sampling}) that this approach might produce a low-quality reconstruction and prediction, especially for small values of $r$. So, we take a look at some computationally simple sorting criteria in the following subsections.

\subsubsection{Sorting criterion 1} \label{sec:sc1}
One of the early and straightforward sorting criteria is to order the modes based on their initial amplitudes $\mathbf{b}$ as given by Eq.~\ref{eq:amp}. In particular, we can define the importance index $I_i$ as:
\begin{align}
     I_i^{\mathbf{I}} = |b_i|,
\end{align}
where $|\cdot|$ denotes the magnitude. However, this mechanism does not take into account that for specific systems, there might be modes with large initial amplitudes but rapidly decaying and/or modes with small initial amplitudes but rapidly growing.

\subsubsection{Sorting criterion 2} \label{sec:sc2}
We can define the modal contribution to the FOM dynamics by considering the quickly decaying modes with high initial amplitude as well as those with low initial amplitude, but rapid growth rate, as both cases are significant to be taken into account. Thus, the importance index $I_i$ can be written as:
\begin{equation} \label{sc1_1}
    I_i^{\mathbf{II}} = |b_i| \left( e^{\sigma_i} + e^{-\sigma_i} \right),
\end{equation}
where $\sigma_i = \mbox{real} (\alpha_i)$ stands for the associated growth or decay rate of the particular DMD mode. This criterion depends on the initial amplitude $b_i$ and the growth rate $\sigma_i$ of the modes. It was shown that this criterion provides improved results for similar SWE systems on Cartesian coordinates \cite{ahmed2020sampling}.

\subsubsection{Sorting criterion 3} \label{sec:sc3}
\citet{kou2017improved} provide a criterion for evaluating the contribution of DMD modes and sorting them by estimating the influence of each DMD mode over the entire sampling window based on their temporal evolution. The influence of the modes $I_i$ can be evaluated as:
\begin{equation} \label{eq:SC2_1}
    I_i^{\mathbf{III}} = \left( \sum_{j=1}^{m} | b_i \lambda_i^{j-1}| \| \psi_i\|_F^2 \right) \Delta T,
\end{equation}
where $\lambda_i$ are the discrete-time eigenvalues. The subscript $i$ in $\lambda$ refers to the $i^{th}$ mode, while the superscript $(j-1)$ is the power and $m$ is the number of collected snapshots. 

\subsubsection{Sorting criterion 4} \label{sec:sc4}
We also introduce a scheme for sorting the DMD modes based on their initial amplitudes and their growth/decay rates, similar to the one proposed in \cite{kou2017improved}. The general idea is to obtain an estimation of each mode's contribution from the projection of the FOM onto that DMD mode (or basis) over the entire sampling window. In Eq.~\ref{eq:recon_cont}, $e^{\alpha_i t} b_i$ can be thought of as the projection of the FOM onto the $i^{th}$ DMD mode at time $t$. Taking an integral with respect to time of the absolute of this quantity over the entire sampling window ($T$) gives an average contribution of the modes for the given sampling window as follows:
\begin{equation} \label{eq:sc4_1}
    I_i^{\mathbf{IV}} = \frac{1}{T} \int_0^{T} |e^{\alpha_i t} b_i| dt,
\end{equation}
where $T=(m-1) \Delta T$ is the total sampling time and $m$ is the number of snapshots available. Then, Eq.~\ref{eq:sc4_1} can be approximated as follows:
\begin{equation}
    I_i^{\mathbf{IV}} = |b_i|  \frac{e^{\sigma_i T} - 1 }{\sigma_i T}.
\end{equation}

It is interesting to note here that the performance of a mode selection criterion is dependent on the dynamics of the system it is applied to as well. All DMD modes can be assigned a value $I_i$ based on their contribution to the FOM and sorted in decreasing order of $I_i$. Then, the first $r$ modes are selected based on their importance index $I_i$. We also highlight that physically conserving sorting and selection criteria might be tailored for the specific problem in hand. For example, in the case of SWEs, conservation of integral invariants  (e.g., total mass, total energy and potential enstrophy) can be enforced as hard constraints for the optimization and selection of the most important modes \cite{bistrian2017method}. In other words, if we know something about the system being analyzed (e.g., constitutive laws), we can enforce the sorting criteria to improve the DMD results.

\section{Sketchy Dynamic Mode Decomposition} \label{sec:sketch}
As we discussed before, one of the main bottlenecks of the computational pipeline of the stable DMD implementation is the size of the given data matrix. In particular, for high-dimensional systems, the size of the input data matrix $\mathbf{X}$ can be very large and loading $\mathbf{X}_1$ to perform SVD becomes computationally expensive. Streaming algorithms can be adopted to mitigate this burden by passing the snapshots one-by-one and performing SVD on increments of the data \cite{hemati2014dynamic,pendergrass2016streaming}. Another approach, that we consider in the present study, is to utilize sketching as a tool to reduce the size of the processed data matrix. In particular, we seek a low-order embedding of the original data matrix represented by a smaller matrix, called the sketch that captures the main information of the original one. The success of the sketching method relies on the assumption that big data matrices are often low-rank \cite{udell2019big}, with an exponential decay of the underlying singular values. Indeed, this is the main reason why model order reduction has witnessed substantial success in many applications. 

Sketching-based algorithms exploit informed projections to transform the original data matrix to a more compact one. The expensive computations are thus performed onto the latter, after which a post-processing takes place to map (and maybe correct) the outputs to the original space. Randomized projections have been especially effective for this purpose. They can efficiently and accurately extract the spatiotemporal coherent structures from high-dimensional data, \emph{with high probability}. We consider three variants of sketching-based DMD, to reduce the computational costs of different portions of the deterministic DMD algorithm. 

\subsection{Sketching the range of $\mathbf{X}_1$} \label{sec:bist}
\citet{bistrian2017randomized} introduced a randomized DMD framework that aims at mitigating the cost of the SVD of the data matrix $\mathbf{X}_1 \in \mathbb{R}^{n\times(m-1)}$. In particular, a near-optimal basis with a target rank $k$ is defined using random projections such that it captures the range of $\mathbf{X}_1$ and provides a smaller sketch matrix $\mathbf{B}_1 \in \mathbb{R}^{k\times m-1}$. This results in increased efficiency in terms of computational memory and/or time for later deterministic steps of model reduction.

In randomized DMD, a near-optimal orthonormal basis $\mathbf{Q}$ is estimated from the full order dense data matrix $\mathbf{X}_1 = \{ \mathbf{x^{(1)}}, \mathbf{x^{(2)}}, \dots,\mathbf{x^{(m-1)}} \} \in \mathbb{R}^{n \times (m-1)}$ such that $\mathbf{X}_1 \approx \mathbf{Q} \mathbf{Q}^* \mathbf{X}_1$. In order to do this, a randomized projection of the input data matrix $\mathbf{X}_1$ is first performed as follows:
\begin{equation}
    \mathbf{Y}_1 = \mathbf{X}_1 \boldsymbol{\Omega}_1,
\end{equation}
where $\boldsymbol{\Omega}_1 \in \mathbb{R}^{(m-1) \times k}$ is a random matrix drawn from Gaussian distribution which incorporates the randomized concept while $\mathbf{Y}_1 \in \mathbb{R}^{n \times k}$ is a summary of the action of $\mathbf{X}_1$. We highlight that the computational cost of this matrix-matrix multiplication can be reduced by exploiting structured randomized matrices. Here, $k=r+s$ represents the target rank where $s$ is defined as the oversampling factor that helps to obtain an improved basis. \citet{bistrian2017randomized} suggested using an oversampling factor of $s=r$, which we adopt in the present study.

Next, a QR decomposition of $\mathbf{Y}_1=\mathbf{Q}_1\mathbf{R}_1$ is performed where $\mathbf{Q}_1 \in \mathbb{R}^{n \times k}$ and $\mathbf{R}_1 \in \mathbb{R}^{k \times k}$. Indeed, we only need $\mathbf{Q}$ for our computations, and can safely discard $\mathbf{R}_1$. The full data $\mathbf{X}_1$ is then projected onto the basis $\mathbf{Q}_1$ to obtain a lower dimension matrix $\mathbf{B}_1 \in \mathbb{R}^{k \times (m-1)}$ as follows: 
\begin{equation}
    \mathbf{B}_1 = \mathbf{Q}_1^{*} \mathbf{X}_1,
\end{equation}
where $\mathbf{Q}^{*}$ is the conjugate transpose of $\mathbf{Q}$. Thus, SVD can be performed on $\mathbf{B}_1$ instead of $\mathbf{X}_1$ as $\mathbf{B}_1 = \tilde{\mathbf{U}} \tilde{\boldsymbol{\Sigma}} \tilde{\mathbf{V}}^*$. The singular values and vectors of $\mathbf{X}_1$ can be recovered as follows:
\begin{align}
    \mathbf{U} &= \mathbf{Q}_1 \tilde{\mathbf{U}}, \\
    \boldsymbol{\Sigma} &= \tilde{\boldsymbol{\Sigma}}, \\
    \mathbf{V} &=  \tilde{\mathbf{V}}.
\end{align}
The algorithmic steps for the randomized DMD based on sketching the range of $\mathbf{X}_1$ are summarized in Algorithm~\ref{alg:bist}

\begin{algorithm}[H]
\caption{Randomized Dynamic Mode Decomposition by Sketching the Range of $\mathbf{X}_1$}
\label{alg:bist}
\begin{algorithmic}[1]
\STATE The matrix $\mathbf{X}$ is split into two matrices $\mathbf{X}_1 = \{ \mathbf{x}^{(1)}, \mathbf{x}^{(2)}, \dots, \mathbf{x}^{(m-1)} \}$ and $\mathbf{X}_2 = \{ \mathbf{x}^{(2)}, \mathbf{x}^{(3)}, \dots, \mathbf{x}^{(m)} \}$. \\
\vspace{10pt}
\STATE Draw a random matrix $\boldsymbol{\Omega}_1 \in \mathbb{R}^{(m-1) \times k}$ from Gaussian distribution and perform the randomized projection of $\mathbf{X}_1$ \vspace{-10pt}
\begin{align*}
     \mathbf{Y}_1 = \mathbf{X}_1 \boldsymbol{\Omega}_1
\end{align*}

\STATE Perform QR decomposition as $\mathbf{Y}_1=\mathbf{Q}_1\mathbf{R}_1$ to obtain a near-optimal basis $\mathbf{Q}_1$ for $\mathbf{X}_1$ and discard $\mathbf{R}_1$.
\vspace{10pt}
\STATE A sketch $\mathbf{B}_1$ of $\mathbf{X}_1$ is obtained as  \vspace{-10pt}
\begin{align*}
        \mathbf{B}_1 = \mathbf{Q}_1^{*} \mathbf{X}_1
\end{align*}

\STATE Perform SVD on $\mathbf{B}_1$ \vspace{-10pt}
\begin{align*}
    \tilde{\mathbf{U}} \tilde{\boldsymbol{\Sigma}} \tilde{\mathbf{V}}^* = \text{svd} {\left( \mathbf{B}_1 \right)}
\end{align*}

\STATE Recover SVD of $\mathbf{X}_1$ \vspace{-10pt}
\begin{align*}
    \mathbf{U} &= \mathbf{Q}_1 \tilde{\mathbf{U}} \\
    \boldsymbol{\Sigma} &= \tilde{\boldsymbol{\Sigma}} \\
    \mathbf{V} &=  \tilde{\mathbf{V}}
\end{align*}

\STATE Rank truncation [to reduce noise] \vspace{-10pt}
\begin{align*}
    \mathbf{U}_\mathbf{R} &= \mathbf{U}(:,1:R) \\
    \mathbf{V}_\mathbf{R} &= \mathbf{V}(:,1:R) \\
    \boldsymbol{\Sigma}_\mathbf{R} &= \boldsymbol{\Sigma}(1:R,1:R)
\end{align*}

\STATE Low-rank dynamics \vspace{-10pt}
\begin{align*}
    \Tilde{\mathbf{A}} = \mathbf{U^*_R} \mathbf{X}_2 \mathbf{V}^{}_\mathbf{R} \boldsymbol{\Sigma}^{-1}_\mathbf{R}
\end{align*}

\STATE Eigenvalue decomposition \vspace{-10pt}
\begin{align*}
    \left[\mathbf{W}, \boldsymbol{\Lambda} \right] = \text{eig} {( \Tilde{\mathbf{A}} )}
\end{align*}

\STATE Compute DMD modes and spectrum \vspace{-10pt}
\begin{align*}
\boldsymbol{\Psi} &= \mathbf{U}_\mathbf{R} \mathbf{W} \quad \text{or} \quad \boldsymbol{\Psi} = \mathbf{X}_2 \mathbf{V}^{}_\mathbf{R} \boldsymbol{\Sigma}^{-1}_\mathbf{R} \mathbf{W} \\
\lambda_i &= \{ \text{diag} (\boldsymbol{\Lambda}) \} \\
\alpha_i &= \text{ln} {(\lambda_i)}/\Delta T
\end{align*}
\end{algorithmic}
\end{algorithm}

\subsection{Sketching the range of $\mathbf{X}$} \label{sec:erich}
In Section~\ref{sec:bist}, the deterministic SVD of $\mathbf{X}_1$ is bypassed by a randomized projection of $\mathbf{X}_1$ to obtain a sketch that captures its range and efficiently applying SVD onto this sketch. After the SVD of $\mathbf{X}_1$ is recovered, the remaining steps are performed in a high-dimensional space. Alternatively, \citet{erichson2019randomized} developed a randomized DMD algorithm that is based on sketching the range of $\mathbf{X}$ and performing \emph{all steps} in the DMD procedure in a low-order space, while the DMD of the original system is recovered at the very end.
Therefore, a random projection is defined to capture the range of $\mathbf{X}$ as follows:
\begin{equation}
    \mathbf{Y} = \mathbf{X} \boldsymbol{\Omega},
\end{equation}
where $\boldsymbol{\Omega} \in \mathbb{R}^{m \times k}$ is a random matrix and $\mathbf{Y} \in \mathbb{R}^{n \times k}$ is the summary of the action of $\mathbf{X}$. Next, a QR decomposition of $\mathbf{Y}=\mathbf{Q}\mathbf{R}$ is performed where $\mathbf{Q} \in \mathbb{R}^{n \times k}$ and $\mathbf{R} \in \mathbb{R}^{k \times k}$, where $\mathbf{Q}$ is considered a near-optimal orthonormal basis such that $\mathbf{X} \approx \mathbf{Q} \mathbf{Q}^* \mathbf{X}$. The full data $\mathbf{X}$ is projected onto $\mathbf{Q}$ to obtain a lower dimension matrix $\mathbf{B} \in \mathbb{R}^{k \times m}$ as follows:
\begin{equation}
    \mathbf{B} = \mathbf{Q}^{*} \mathbf{X}.
\end{equation}
$\mathbf{B}$ can be further split into two matrices $\mathbf{B}_1 \in \mathbb{R}^{k \times (m-1)}$ and $\mathbf{B}_2 \in \mathbb{R}^{k \times (m-1)}$ by selecting the first and last $(m-1)$ columns of $\mathbf{B}$, respectively.

The \emph{whole} DMD algorithm can be applied onto the low-order sketches to yield the DMD modes $\tilde{\boldsymbol{\Psi}}$ of the low-order matrix $\mathbf{B}$, and the DMD modes are recovered at the end as follows:
\begin{align}
    \boldsymbol{\Psi} = \mathbf{Q} \tilde{\boldsymbol{\Psi}}
\end{align}
The algorithmic steps for the randomized DMD based on sketching the range of $\mathbf{X}$ are summarized in Algorithm~\ref{alg:erich}

\begin{algorithm}[H]
\caption{Randomized Dynamic Mode Decomposition by Sketching the Range of $\mathbf{X}$}
\label{alg:erich}
\begin{algorithmic}[1]

\STATE Draw a random matrix $\boldsymbol{\Omega} \in \mathbb{R}^{m \times k}$ from Gaussian distribution and perform the randomized projection of $\mathbf{X}$ \vspace{-10pt}
\begin{align*}
     \mathbf{Y} = \mathbf{X} \boldsymbol{\Omega}
\end{align*}

\STATE Perform QR decomposition as $\mathbf{Y}=\mathbf{Q}\mathbf{R}$ to obtain a near-optimal basis $\mathbf{Q}$ for $\mathbf{X}$ and discard $\mathbf{R}$. 
 \vspace{10pt}
\STATE A sketch $\mathbf{B}$ of $\mathbf{X}$ is obtained as  \vspace{-10pt}
\begin{align*}
        \mathbf{B} = \mathbf{Q}^{*} \mathbf{X}
\end{align*}

\STATE The matrix $\mathbf{B}$ is split into two matrices $\mathbf{B}_1 \in \mathbb{R}^{k \times (m-1)}$ and $\mathbf{B}_2 \in \mathbb{R}^{k \times (m-1)}$.
\vspace{10pt}
\STATE Perform SVD on $\mathbf{B}_1$ \vspace{-10pt}
\begin{align*}
    \tilde{\mathbf{U}} \tilde{\boldsymbol{\Sigma}} \tilde{\mathbf{V}}^* = \text{svd} {\left( \mathbf{B}_1 \right)}
\end{align*}

\STATE Rank truncation [to reduce noise] \vspace{-10pt} \vspace{-10pt}
\begin{align*}
    \tilde{\mathbf{U}}_\mathbf{R} &= \tilde{\mathbf{U}}(:,1:R) \\
    \tilde{\mathbf{V}}_\mathbf{R} &= \tilde{\mathbf{V}}(:,1:R) \\
    \tilde{\boldsymbol{\Sigma}}_\mathbf{R} &= \tilde{\boldsymbol{\Sigma}}(1:R,1:R) 
\end{align*}

\STATE Low-rank dynamics \vspace{-10pt}
\begin{align*}
    \Tilde{\mathbf{A}} = \tilde{\mathbf{U}}^*_\mathbf{R} \mathbf{B}_2 \tilde{\mathbf{V}}_\mathbf{R} \tilde{\boldsymbol{\Sigma}}^{-1}_\mathbf{R}
\end{align*}

\STATE Eigenvalue decomposition \vspace{-10pt}
\begin{align*}
    \left[\mathbf{W}, \boldsymbol{\Lambda} \right] = \text{eig} {( \Tilde{\mathbf{A}} )}
\end{align*}

\STATE Compute DMD modes for $\mathbf{B}$ \vspace{-10pt}
\begin{align*}
\tilde{\boldsymbol{\Psi}} &= \tilde{\mathbf{U}}_R \mathbf{W} \quad \text{or} \quad \tilde{\boldsymbol{\Psi}} = \mathbf{B}_2 \tilde{\mathbf{V}}^{}_\mathbf{R} \tilde{\boldsymbol{\Sigma}}^{-1}_\mathbf{R} \mathbf{W}
\end{align*}

\STATE Recover DMD modes of $\mathbf{X}$ \vspace{-10pt}
\begin{align*}
    \boldsymbol{\Psi} &= \mathbf{Q} \tilde{\boldsymbol{\Psi}}
\end{align*}

\STATE Compute DMD spectrum \vspace{-10pt}
\begin{align*}
\lambda_i &= \{ \text{diag} (\boldsymbol{\Lambda}) \} \\
\alpha_i &= \text{ln} {(\lambda_i)}/\Delta T
\end{align*}
\end{algorithmic}
\end{algorithm} 

\clearpage
\subsection{Sketching the range and corange of $\mathbf{X}_1$} \label{sec:sketch1}
In this study, we further develop an efficient DMD implementation based on a sketching approach that captures the range and corange of the data matrix $\mathbf{X}_1 \in \mathbb{R}^{n\times(m-1)}$, resulting in an even smaller sketch than the range sketches \cite{yurtsever2017sketchy,tropp2019streaming}. We first describe the sketching operators parametrized by a ``range'' parameter $k$ and a ``core'' parameter $p$ that satisfy $r\le k \le p \le \min{(n,m-1)}$, where the parameter $k$ determines the maximum rank of an approximation.

Now, we independently draw and fix four randomized linear reduction maps (often called test matrices) as follows:
\begin{align*}
    \boldsymbol{\Omega}_1 \in \mathbb{R}^{(m-1)\times k} \quad &\text{and} \quad \boldsymbol{\Gamma}_1 \in \mathbb{R}^{k\times n}   \\
   \boldsymbol{\Theta}_1 \in \mathbb{R}^{(m-1)\times p}  \quad &\text{and} \quad \boldsymbol{\Phi}_1 \in \mathbb{R}^{p\times n}.
\end{align*}
Then, we define three matrices comprising our sketch as follows:
\begin{align}
    \mathbf{F}_1 &= \mathbf{X}_1 \boldsymbol{\Omega}_1 \in \mathbb{R}^{n\times k},  \\
    \mathbf{G}_1 &= \boldsymbol{\Gamma}_1 \mathbf{X}_1 \in \mathbb{R}^{k\times (m-1)}, \\
    \mathbf{H}_1 &= \boldsymbol{\Phi}_1 \mathbf{X}_1 \boldsymbol{\Theta}_1 \in \mathbb{R}^{p\times p},
\end{align}
where $\mathbf{F}_1$ and $\mathbf{G}_1$ capture the range and corange of $\mathbf{X}_1$, respectively, while $\mathbf{H}_1$ is called the core sketch and describes how $\mathbf{X}_1$ acts between the spaces captured by sketches $\mathbf{F}_1$ and $\mathbf{G}_1$ \cite{rajaram2020randomized}. Now, near-optimal bases for the range and corange of $\mathbf{X}_1$ are computed by the thin QR factorization of $\mathbf{F}_1$ and $\mathbf{G}_1$ as follows:
\begin{align}
    \mathbf{F}_1 &= \mathbf{Q}_1 \mathbf{R}_1, \\
    \mathbf{G}^*_1 &= \mathbf{P}_1 \mathbf{T}_1,
\end{align}
where we can again discard the triangular matrices $\mathbf{T}_1$ and $\mathbf{R}_1$. The third sketch $\mathbf{H}_1$ is used to compute the core approximation $\mathbf{C}_1 \in \mathbb{R}^{k\times k}$  of $\mathbf{X}_1$ as follows:
\begin{align}
    \mathbf{C}_1 = (\boldsymbol{\Phi}_1 \mathbf{Q}_1)^{\dagger} \mathbf{H}_1 ( \mathbf{P}_1^* \boldsymbol{\Theta}_1) ^{\dagger},
\end{align}
which can be implemented efficiency by solving a family of least-squares problems. The original data matrix is now related to the core sketch by the following relation:
\begin{align}
    \mathbf{X}_1 = \mathbf{Q}_1 \mathbf{C}_1 \mathbf{P}^*_1.
\end{align}

Note that $\mathbf{C}_1 \in \mathbb{R}^{k\times k}$ is a square matrix with small size $k$, compared to $\mathbf{B}_1$ and $\mathbf{B}$ in Section~\ref{sec:bist}-\ref{sec:erich}. Thus, we use the core sketch to compute the SVD of $\mathbf{C}_1 = \tilde{\mathbf{U}} \tilde{\boldsymbol{\Sigma}} \mathbf{V}^*$ efficiently and the singular values and vectors of $\mathbf{X}_1$ can be recovered as follows:
\begin{align}
    \mathbf{U} &= \mathbf{Q}_1 \tilde{\mathbf{U}}, \\
    \boldsymbol{\Sigma} &= \tilde{\boldsymbol{\Sigma}}, \\
    \mathbf{V} &=  \mathbf{P}_1 \tilde{\mathbf{V}}.
\end{align}
The algorithmic steps for the randomized DMD based on sketching the range of $\mathbf{X}_1$ are summarized in Algorithm~\ref{alg:sketch1}

\begin{algorithm}[H]
\caption{Randomized Dynamic Mode Decomposition by Sketching the Range and Corange of $\mathbf{X}_1$}
\label{alg:sketch1}
\begin{algorithmic}[1]
\STATE The matrix $\mathbf{X}$ is split into two matrices $\mathbf{X}_1 = \{ \mathbf{x}^{(1)}, \mathbf{x}^{(2)}, \dots, \mathbf{x}^{(m-1)} \}$ and $\mathbf{X}_2 = \{ \mathbf{x}^{(2)}, \mathbf{x}^{(3)}, \dots, \mathbf{x}^{(m)} \}$. \\
\vspace{10pt}
\STATE Draw random matrices $\boldsymbol{\Omega}_1 \in \mathbb{R}^{(m-1)\times k}$, $\boldsymbol{\Gamma}_1 \in \mathbb{R}^{k\times n}$, $\boldsymbol{\Theta}_1 \in \mathbb{R}^{(m-1)\times s}$, and $\boldsymbol{\Phi}_1 \in \mathbb{R}^{s\times n}$ independently from Gaussian distribution and perform the randomized projection of $\mathbf{X}_1$ \vspace{-10pt}
\begin{align*}
    \mathbf{F}_1 &= \mathbf{X}_1 \boldsymbol{\Omega}_1 \in \mathbb{R}^{n\times k},  \\
    \mathbf{G}_1 &= \boldsymbol{\Gamma}_1 \mathbf{X}_1 \in \mathbb{R}^{k\times (m-1)}, \\
    \mathbf{H}_1 &= \boldsymbol{\Phi}_1 \mathbf{X}_1 \boldsymbol{\Theta}_1 \in \mathbb{R}^{s\times s},
\end{align*}

\STATE Perform QR decomposition of $\mathbf{F}_1$ and $\mathbf{G}_1$ as $\mathbf{F}_1 = \mathbf{Q}_1 \mathbf{R}_1$ and $\mathbf{G}^*_1 = \mathbf{P}_1 \mathbf{T}_1$ to obtain a near-optimal basis $\mathbf{Q}_1$ and $\mathbf{P}_1$ for the range and corange of $\mathbf{X}_1$, respectively. Discard $\mathbf{T}_1$ and $\mathbf{R}_1$.
\vspace{10pt}
\STATE A core sketch $\mathbf{C}_1 \in \mathbb{R}^{k\times k}$  of $\mathbf{X}_1$ is obtained as follows, \vspace{-10pt}
\begin{align*}
\mathbf{C}_1 = (\boldsymbol{\Phi}_1 \mathbf{Q}_1)^{\dagger} \mathbf{H}_1 ( \mathbf{P}_1^* \boldsymbol{\Theta}_1) ^{\dagger}
\end{align*}

\STATE Perform SVD on $\mathbf{C}_1$ \vspace{-10pt}
\begin{align*}
    \tilde{\mathbf{U}} \tilde{\boldsymbol{\Sigma}} \tilde{\mathbf{V}}^* = \text{svd} {\left( \mathbf{C}_1 \right)}
\end{align*}

\STATE Recover SVD of $\mathbf{X}_1$ \vspace{-10pt}
\begin{align*}
    \mathbf{U} &= \mathbf{Q}_1 \tilde{\mathbf{U}} \\
    \boldsymbol{\Sigma} &= \tilde{\boldsymbol{\Sigma}} \\
    \mathbf{V} &=  \mathbf{P}_1 \tilde{\mathbf{V}}
\end{align*}

\STATE Rank truncation [to reduce noise] \vspace{-10pt}
\begin{align*}
    \mathbf{U}_\mathbf{R} &= \mathbf{U}(:,1:R) \\
    \mathbf{V}_\mathbf{R} &= \mathbf{V}(:,1:R) \\
    \boldsymbol{\Sigma}_\mathbf{R} &= \boldsymbol{\Sigma}(1:R,1:R)
\end{align*}

\STATE Low-rank dynamics \vspace{-10pt}
\begin{align*}
    \Tilde{\mathbf{A}} = \mathbf{U^*_R} \mathbf{X}_2 \mathbf{V}^{}_\mathbf{R} \boldsymbol{\Sigma}^{-1}_\mathbf{R}
\end{align*}

\STATE Eigenvalue decomposition \vspace{-10pt}
\begin{align*}
    \left[\mathbf{W}, \boldsymbol{\Lambda} \right] = \text{eig} {( \Tilde{\mathbf{A}} )}
\end{align*}

\STATE Compute DMD modes and spectrum \vspace{-10pt}
\begin{align*}
\boldsymbol{\Psi} &= \mathbf{U}_\mathbf{R} \mathbf{W} \quad \text{or} \quad \boldsymbol{\Psi} = \mathbf{X}_2 \mathbf{V}^{}_\mathbf{R} \boldsymbol{\Sigma}^{-1}_\mathbf{R} \mathbf{W} \\
\lambda_i &= \{ \text{diag} (\boldsymbol{\Lambda}) \} \\
\alpha_i &= \text{ln} {(\lambda_i)}/\Delta T
\end{align*}
\end{algorithmic}
\end{algorithm} 
\vspace{-20pt}
\section{Results and Discussions} \label{sec:res}

We consider the numerical solutions of the SWEs on the sphere to generate our data sets. For the initial and boundary conditions defined in Section~\ref{sec:math}, we carry out the simulation a time period of $6$ days ($144$ hrs). The solution of the SWEs gives the velocity field ($u_{\phi}$ and $u_{\theta}$) directly, while the vorticity field ($\omega(\phi,\theta)$) is numerically derived by Eq.~\ref{eq:vorticity2}. Snapshots of the vorticity field data from the numerical solution is utilized for the demonstration of the DMD-ROMs. For the upcoming DMD analysis, snapshots from the end of the $3^{rd}$ day ($t=72$ hr) to the end of the $6^{th}$ day ($t=144$ hr), i.e., a period of $3$ days, are considered for DMD construction. This time frame is henceforth referred to as the sampling window for the DMD analysis. As mentioned in Section \ref{sec:numerics}, the time interval between every consecutive snapshots is $15$ minutes, which shall be referred to as the sampling frequency.

\begin{figure}[H]
    \centering
    \includegraphics[width=0.7\linewidth]{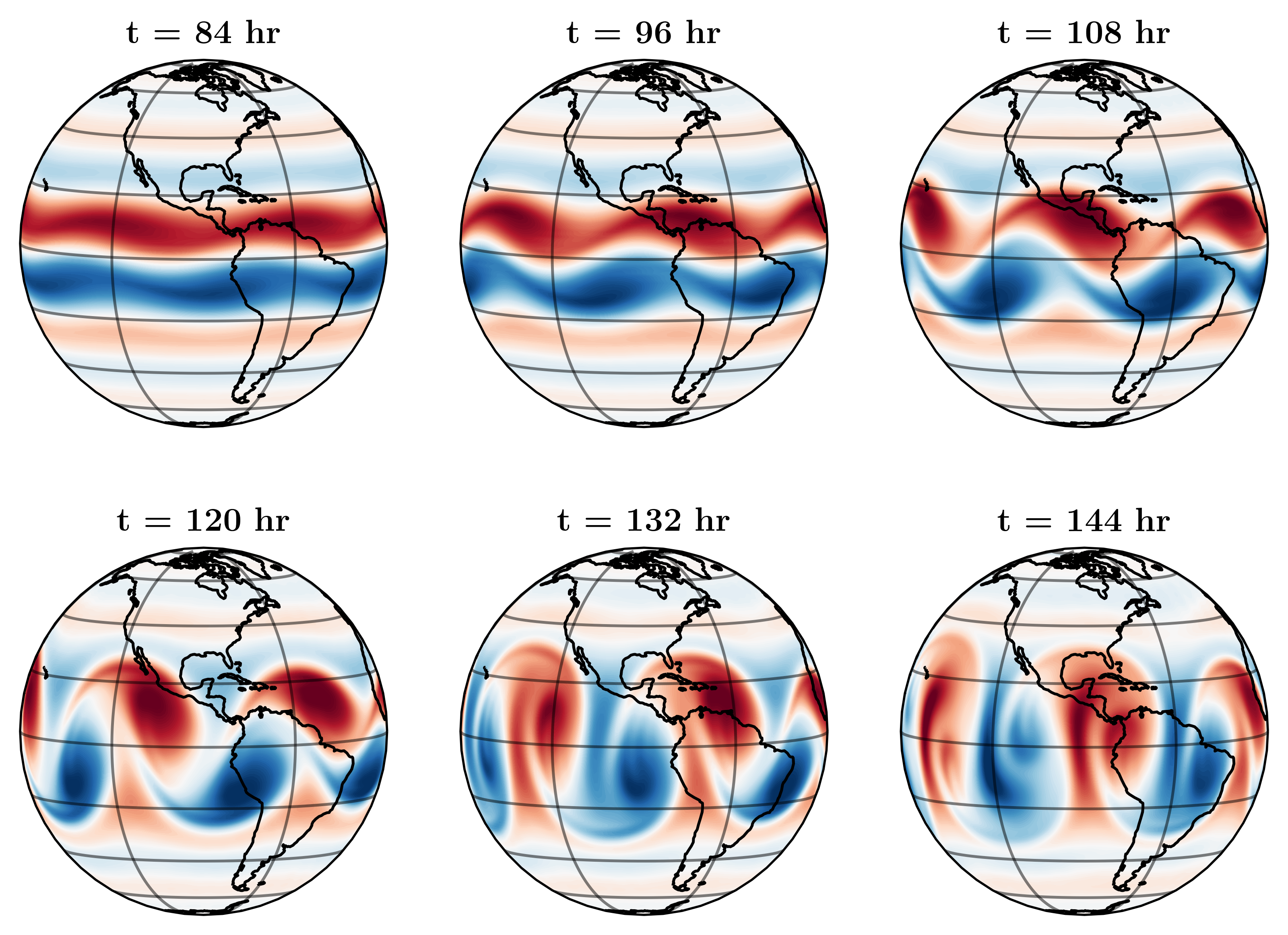}
    \caption{Evolution of the vorticity over the Earth's surface as predicted by the SWEs from $t=84$ hr to $t=144$ hr.} \vspace{-20pt}
    \label{fig:fom}
\end{figure}

The results of numerical simulation for vorticity within the sampling window are presented in Figure \ref{fig:fom}. We utilize the Cartopy package \cite{Cartopy} for data processing and visualization. Figure \ref{fig:fom} reveals the evolution of the vorticity field after an initial period of three days is passed due to the instability in the viscous shear layer arising from the random forcing introduced to the equilibrium condition, as presented in Eq.~\ref{eq:IC1.1} of Section~\ref{sec:math}. 

\subsection{Effect of mode selection}
\begin{figure}[H]
    \centering
    \includegraphics[width=0.7\linewidth]{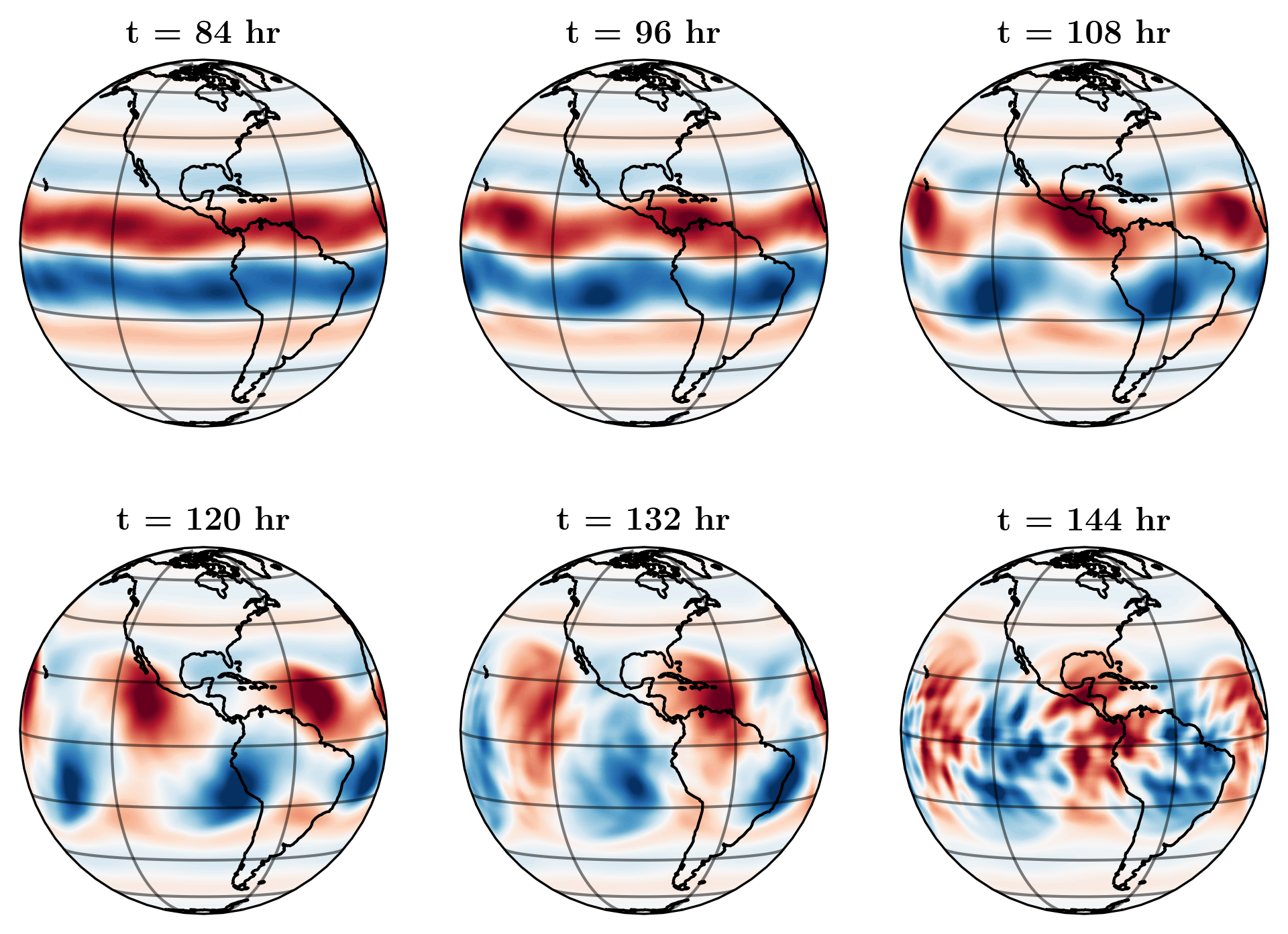}
    \caption{Deterministic DMD predictions with early truncation using $r=20$ modes.}
    \label{fig:det0_20}
\end{figure}

\begin{figure}[H]
    \centering
    \includegraphics[width=0.7\linewidth]{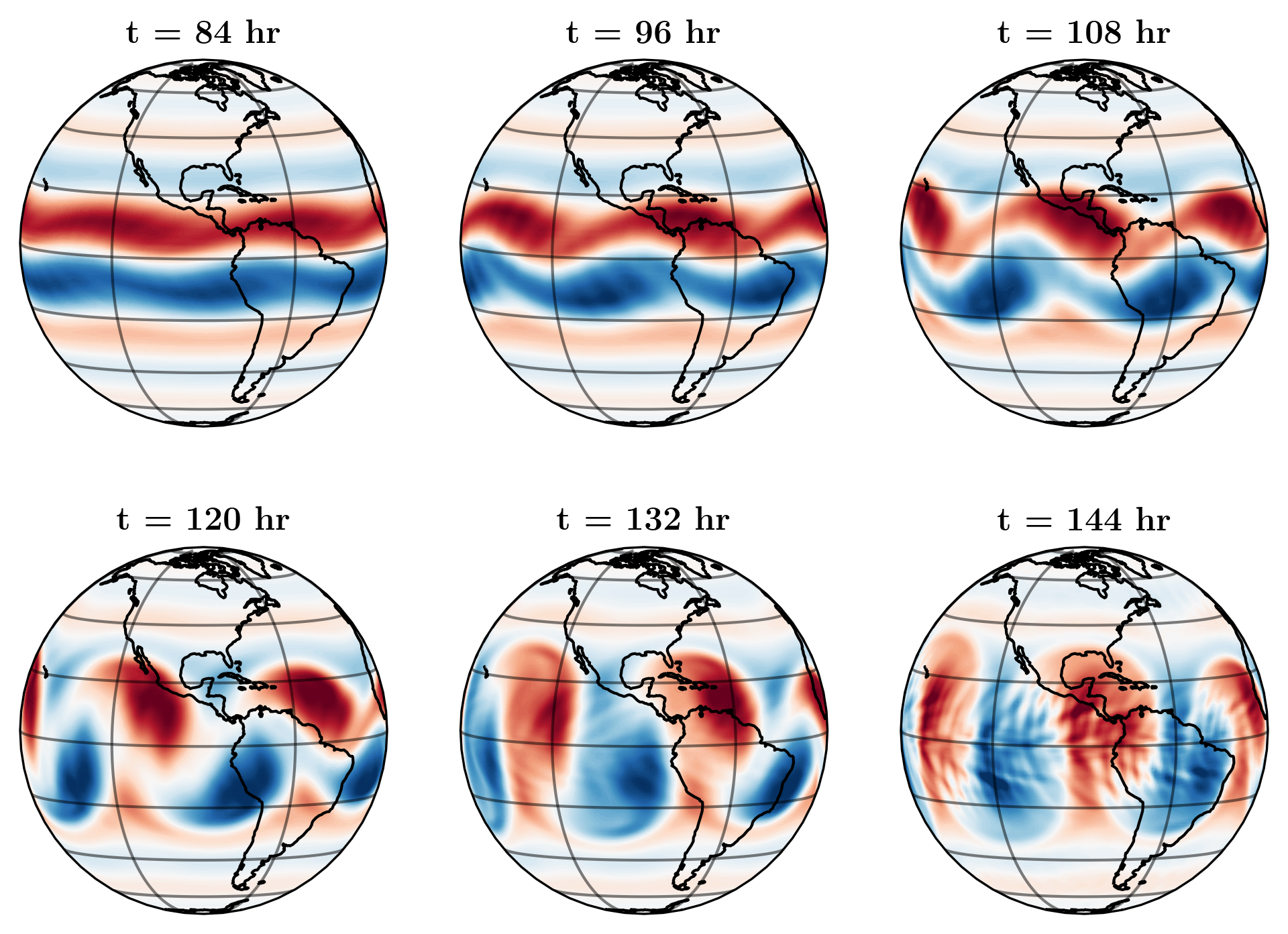}
    \caption{Deterministic DMD predictions with early truncation using $r=40$ modes.} 
    \label{fig:det0_40}
\end{figure}

For the DMD computations, we investigate the effect of the sorting criteria onto the performance of DMD predictions. We first show the performance of the deterministic DMD approach with several selection criteria. We show the prediction of deterministic DMD with early truncation using $r=20$ modes and $r=40$ modes in Figure~\ref{fig:det0_20} and Figure~\ref{fig:det0_40}, respectively. We can observe that an early truncation yields inaccurate results due to the under-estimation of the data matrix $\mathbf{X}_1$, especially with small number of modes.

\begin{figure}[H]
    \centering
    \includegraphics[width=0.7\linewidth]{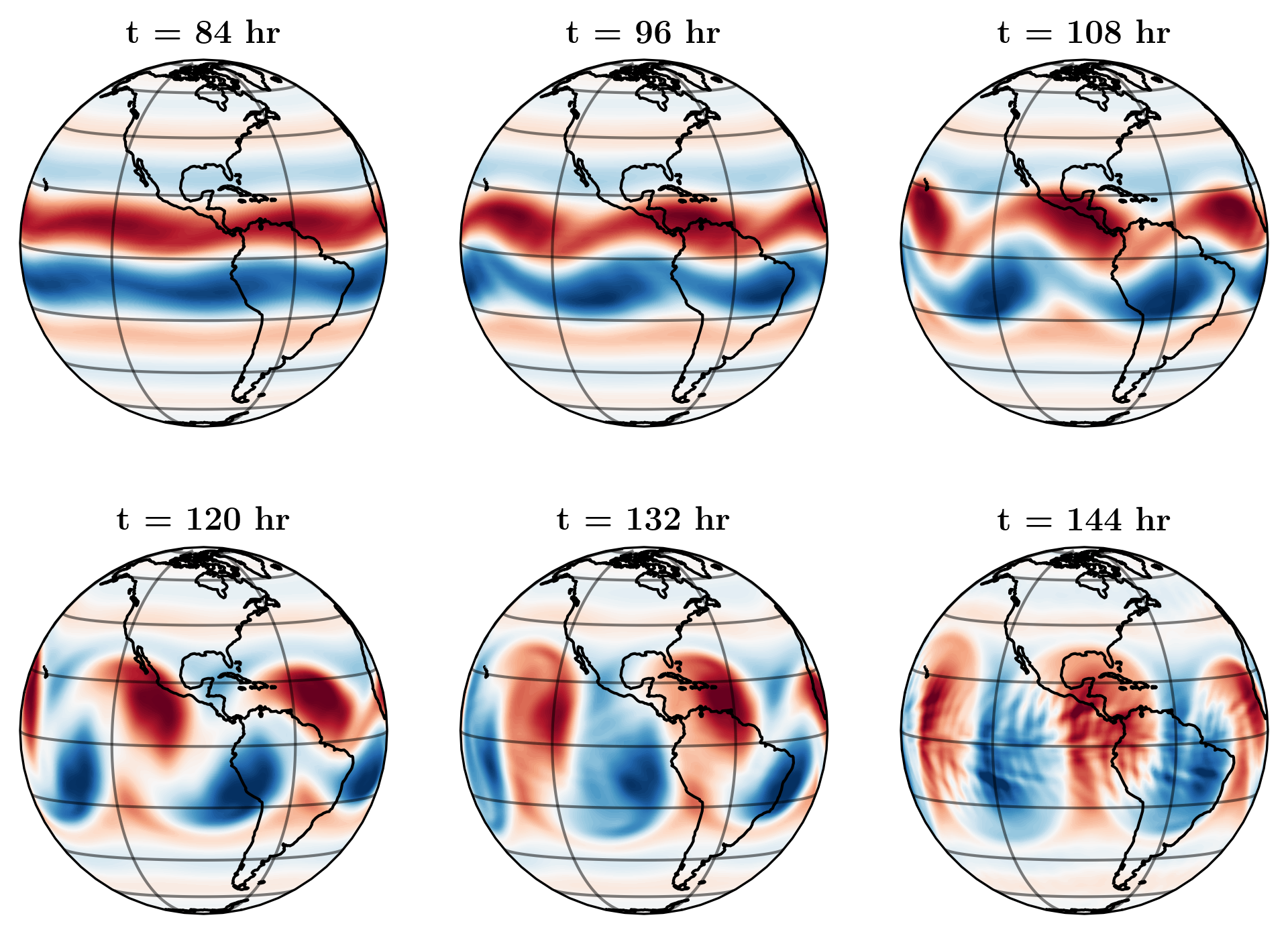}
    \caption{Deterministic DMD predictions with sorting criterion \#4 n using $r=20$ modes.} 
    \label{fig:det4_20}
\end{figure}

On the other hand, with adopting a more \emph{intrusive} selection criterion based on the system's dynamics, we can obtain better results. For example, in Figure~\ref{fig:det4_20}, we demonstrate the predictive capability of deterministic DMD while adopting the sorting criterion \#4. We observe that improved performance is obtained, even comparable to predictions using $40$ modes with early truncation.

\begin{figure}[H]
    \centering
    \includegraphics[width=0.9\linewidth]{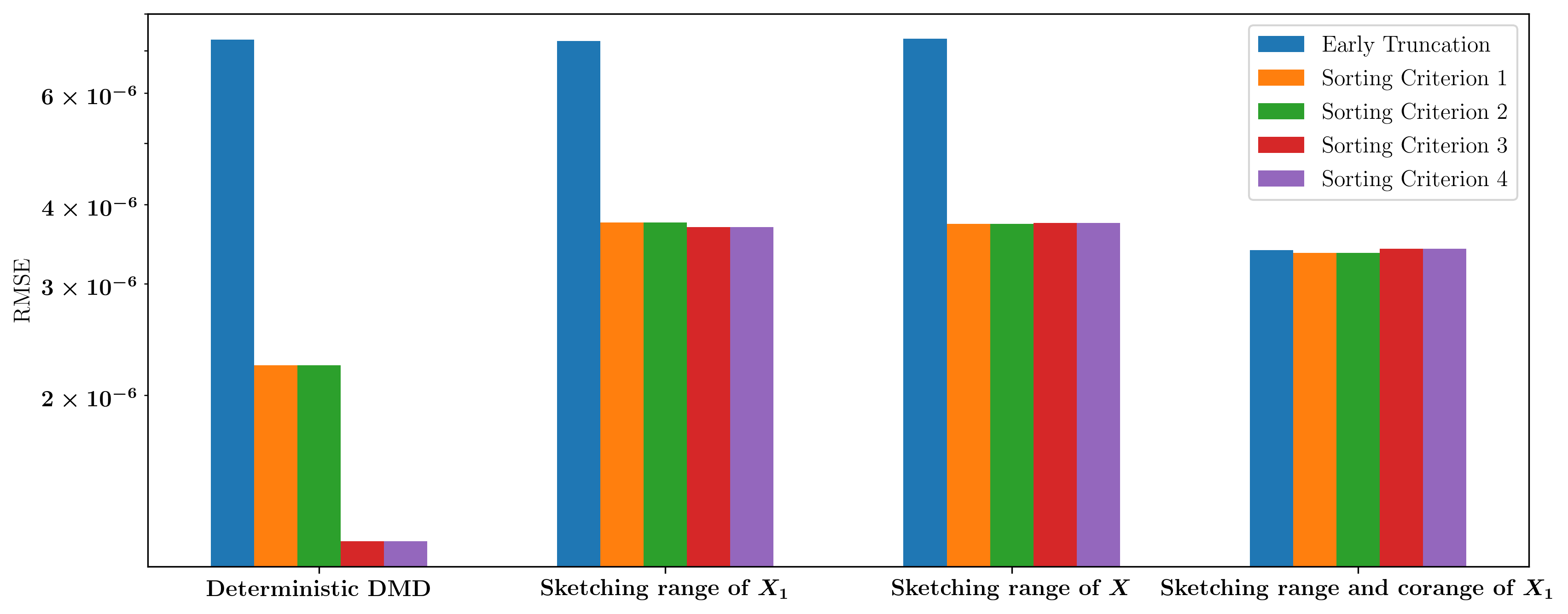} 
    \caption{RMSE of DMD construction with respect to the FOM solution.} 
    \label{fig:sort}
\end{figure}

\clearpage
In order to characterize the results of DMD-ROMs in a quantitative way, we look at the root mean squared error (RMSE). The overall RMSE is mathematically defined as follows
\begin{equation} \label{eq:RMSE}
    \text{RMSE}(t_k) = \sqrt{ \frac{1}{N_{s}} \frac{1}{N_{\phi}} \frac{1}{N_{\theta}} \sum_{i=1}^{N_{\phi}} \sum_{j=1}^{N_{\theta}} \sum_{k=1}^{N_{s}}\left( (\Tilde{\omega}^{(k)}_{i,j})_{FOM} - (\Tilde{\omega}^{(k)}_{i,j})_{ROM} \right)^2}
\end{equation}
where $(\Tilde{\omega}^{k}_{i,j})_{FOM}$ and $(\Tilde{\omega}^{k}_{i,j})_{ROM}$ are the vorticity field data from the full order model and reduced order model by DMD reconstruction, respectively at time $t_k$ at the $i^{th}$ longitude and $j^{th}$ latitude positions. The RMSE values corresponding to different sorting crietria and different DMD implementation are shown in Figure~\ref{fig:sort} for $r=20$ modes. We can see that early truncation yields the least accuracy, except for the DMD with range and corange sketching provided in Section~\ref{sec:sketch}. Nonetheless, the four criteria perform equivalently well for the sketching-based algorithms, so for the sake of brevity, for the rest of the paper we show the results with adopting the sorting criterion \#4.

\subsection{Sketching-based DMD}

\begin{figure}[H]
    \centering
    \includegraphics[width=0.7\linewidth]{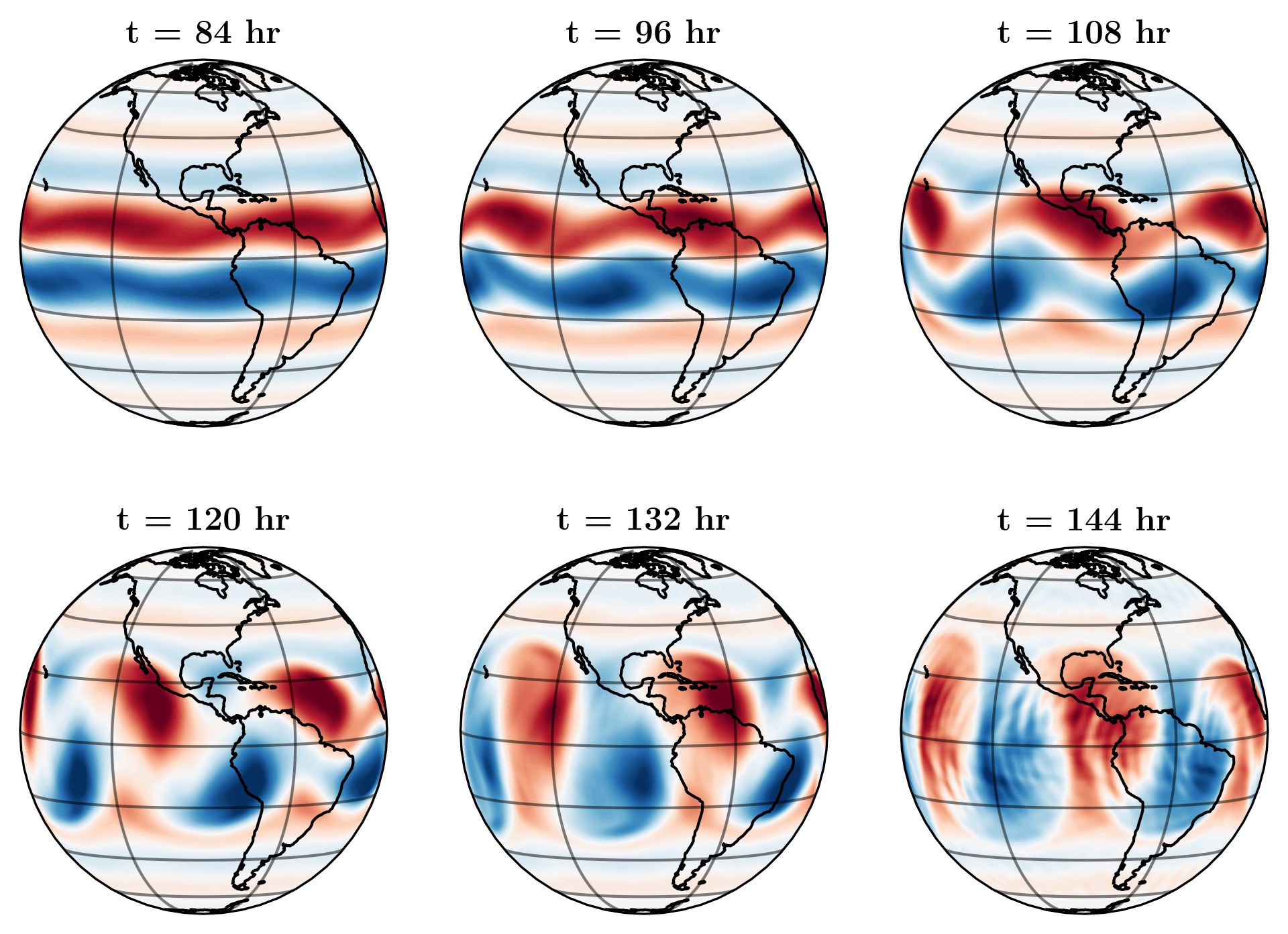}
    \caption{DMD predictions with range of $\mathbf{X}_1$ sketching using $r=20$ modes and sorting criterion \#4.}
    \label{fig:bist4_20}
\end{figure}

We begin by exploring the effect of sketching the range of $\mathbf{X}_1$ to relieve the computational cost of SVD. Figure~\ref{fig:bist4_20} illustrates the performance of the DMD implementation using Algorithm~\ref{alg:bist}, where we can observe that reconstruction quality is comparable to Figure~\ref{fig:det4_20}. Also, Figure~\ref{fig:sort} shows that RMSE is in the same order of magnitude as deterministic DMD. We reiterate that in deterministic DMD, the size of data matrix input to SVD is $n\times (m-1)$ with $n>>m$ while in Algorithm~\ref{alg:bist} the size of the sketch is $k\times(m-1)$ where $k$ is set to $2r$ in the present study.

\begin{figure}[H]
    \centering
    \includegraphics[width=0.7\linewidth]{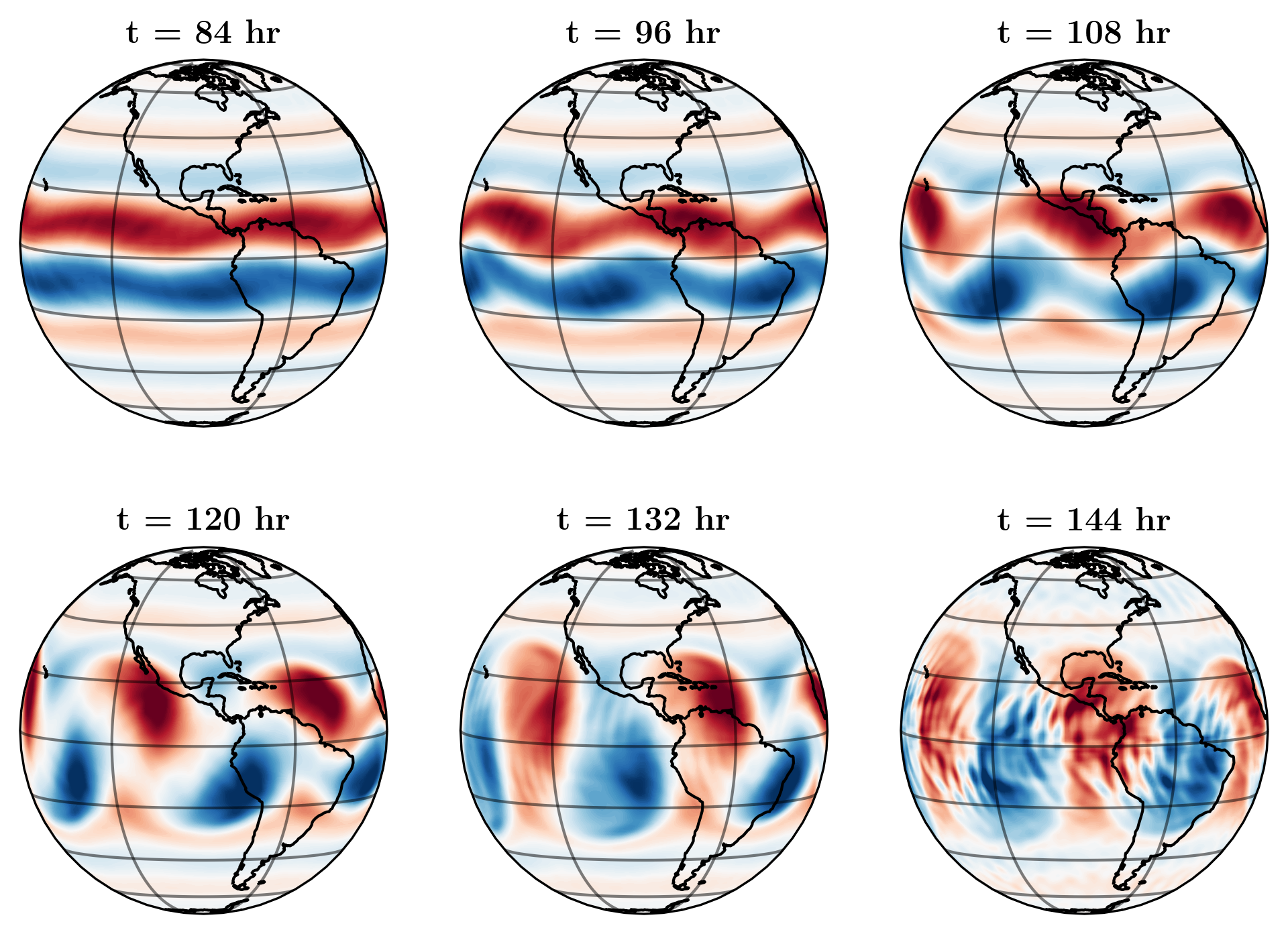}
    \caption{DMD predictions with range of $\mathbf{X}$ sketching using $r=20$ modes and sorting criterion \#4.}
    \label{fig:erich4_20}
\end{figure}

In Figure~\ref{fig:erich4_20}, the outputs of implementing Algorithm~\ref{alg:erich} are shown. In this implementation, a sketch of the range of the data matrix $\mathbf{X}$ is derived and the whole DMD process is performed onto this small-size sketch resulting in additional computational savings, compared to Algorithm~\ref{alg:bist}.

\begin{figure}[H]
    \centering
    \includegraphics[width=0.75\linewidth]{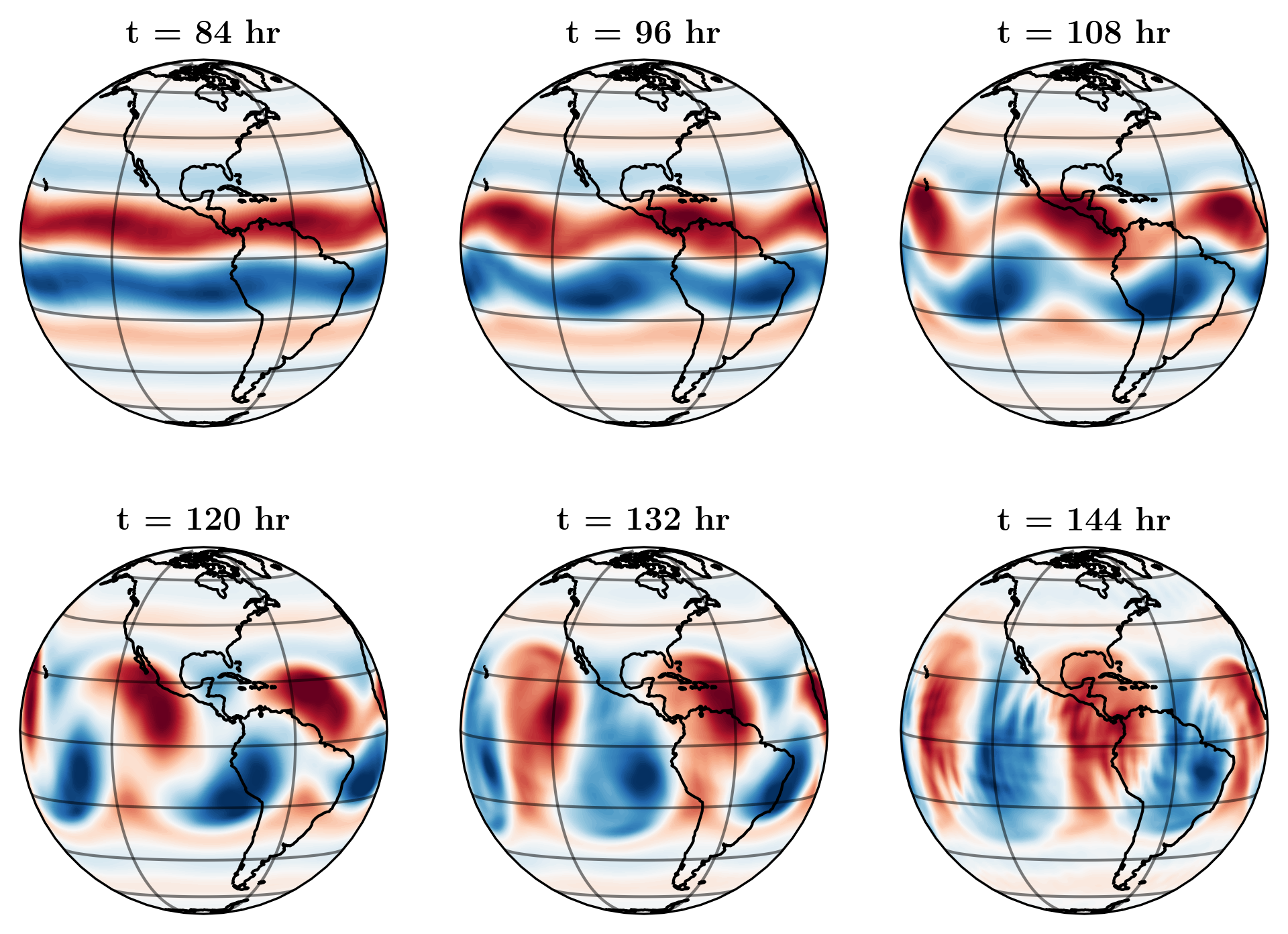}
    \caption{DMD predictions with sketching the range and corange of $\mathbf{X}_1$ using $r=20$ modes and sorting criterion \#4.}
    \label{fig:sketch4_20}
\end{figure}

Finally, the results of DMD based on sketching the range and corange of $\mathbf{X}_1$ with $r=20$ modes are provided in Figure~\ref{fig:sketch4_20}, where we can notice its predictive capabilities compared to other approaches. In this method, the SVD is implemented on a core sketch with a size of $k \times k$ which is usually much smaller than those utilized in Algorithm~\ref{alg:bist} and Algorithm~\ref{alg:erich}. Nonetheless, we highlight that the performance of this approach is quite sensitive to the values of $k$ and $p$. \citet{tropp2019streaming} provide some guidelines on setting these \emph{hyperparameters} given the memory constraints with error bounds guarantees.

\section{Conclusions} \label{sec:conc}
In the present study, we investigate the feasibility of sketching-based algorithms for the efficient implementation of stable singular value decomposition (SVD)-based dynamic mode decomposition (DMD). In particular, we consider three variants based on random projections to sketch the range and corange of the input data matrix. This enables us to bypass expensive portions of the DMD procedure and implement them in a smaller space. Numerical experiments using the spherical shallow water equations are exploited for assessing the applicability of the sketching-based DMD algorithms. We demonstrate that a core sketch with a range and corange projections of the original data matrix can be effectively adopted for the DMD computations. This core sketch is a square matrix with maximum rank of $k$, that can be tuned for specific memory constraints. Moreover, we show that an early truncation of the singular values of the core matrix can yield good results, mitigating the need for further intrusive sorting criterion. Nonetheless, we highlight that the quality of this core sketch is significantly sensitive to the chosen parameters of the sketching algorithm.

\section*{Acknowledgments}
This material is based upon work supported by the U.S. Department of Energy, Office of Science, under the Advanced Scientific Computing Research program (award Number DE-SC0019290) and the National Science Foundation under the Computational Mathematics program (grant DMS-2012255). O.S. gratefully acknowledges their support. \\
Disclaimer: This report was prepared as an account of work sponsored by an agency of the United States Government. Neither the United States Government nor any agency thereof, nor any of their employees, makes any warranty, express or implied, or assumes any legal liability or responsibility for the accuracy, completeness, or usefulness of any information, apparatus, product, or process disclosed, or represents that its use would not infringe privately owned rights. Reference herein to any specific commercial product, process, or service by trade name, trademark, manufacturer, or otherwise does not necessarily constitute or imply its endorsement, recommendation, or favoring by the United States Government or any agency thereof. The views and opinions of authors expressed herein do not necessarily state or reflect those of the United States Government or any agency thereof.
\bibliography{ref}

\end{document}